\documentclass[12pt]{amsart}
\usepackage{amssymb,amsmath,amsthm,mathtools}  
\usepackage{graphicx}
\usepackage{bm,soul}
\usepackage{enumerate}
\usepackage{braket}
\usepackage{amsaddr}
\usepackage{amscd}
\usepackage[dvipsnames]{xcolor}
\usepackage[mathscr]{eucal}
\usepackage{epstopdf}
\usepackage[hidelinks]{hyperref}
\usepackage{subcaption}
\usepackage{cleveref}
\hypersetup{
	colorlinks   = true, 
	urlcolor     = blue, 
	linkcolor    = blue, 
	citecolor   = PineGreen 
}

\numberwithin{equation}{section}

\newcommand\norm[1]{\left\lVert#1\right\rVert}
\newcommand{\tu}{\textup}
\newcommand{\bfa}[1]{\boldsymbol{#1}} 			%

\newcommand{\tm}[1]{{\color{red}{#1}}}

\newcommand{\T}{^{\mbox{\tiny T}}}
\newcommand{\Ts}{^{\,\mbox{\tiny T}}}
\newcommand{\TFC}{Theory of Functional Connections}
\newcommand{\B}[1]{\boldsymbol{#1}}
\newcommand{\cal}[1]{\mathcal{#1}}

\newcommand{\dd}{\; \text{d}}
\DeclareMathOperator*{\minx}{\textbf{min}}

\usepackage[numbers]{natbib}
\usepackage{filecontents}
\usepackage[latin9]{inputenc}
\usepackage[letterpaper]{geometry}
\usepackage{pifont}
\usepackage{float}
\usepackage{color}
\usepackage{adjustbox}
\usepackage{diagbox}
\usepackage{multirow}
\usepackage[makeroom]{cancel}
\usepackage{array}
\usepackage{mathrsfs}
\usepackage{tikz}
\usepackage{comment}
\usepackage{csquotes}

\graphicspath{{./Figures/}}

\usepackage{amsfonts}

\usepackage{color}
\usepackage{adjustbox}
\usepackage{diagbox}
\definecolor{black}{rgb}{0,0,0}

\definecolor{red}{rgb}{1,0,0}

\definecolor{blue}{rgb}{0,0,1}

\numberwithin{equation}{section}

\newcommand{\beq}{\begin{equation}}
\newcommand{\eeq}{\end{equation}}
\newcommand{\beqq}{\begin{equation*}}
\newcommand{\eeqq}{\end{equation*}}
\newcommand{\beqas}{\begin{eqnarray*}}
\newcommand{\eeqas}{\end{eqnarray*}}
\newcommand{\bsp}{\begin{split}}
\newcommand{\eesp}{\end{split}}

\makeatother

\usepackage[english]{babel}

\author[T. M\MakeLowercase{ai and} D. M\MakeLowercase{ortari}]
{\Large T\MakeLowercase{ina} M\MakeLowercase{ai$^{a,b,*}$ and} D\MakeLowercase{aniele} M\MakeLowercase{ortari$^{c}$}} 

\date{\today}

\title[TFC \MakeLowercase{applied to quadratic and nonlinear programming under equality Constraints}]
{\textsf{\LARGE T\MakeLowercase{heory of Functional Connections Applied to Quadratic and Nonlinear Programming under Equality Constraints}$^{\dag}$}\footnote{$^{\dag}$T\MakeLowercase{his paper is an extended version of our published conference paper \cite{mme19}.}}$,^1$\footnote{$^1$A\MakeLowercase{ccepted manuscript by} J\MakeLowercase{ournal of} C\MakeLowercase{omputational and} A\MakeLowercase{pplied} M\MakeLowercase{athematics} (2021).  T\MakeLowercase{he doi of published journal article: \url{https://doi.org/10.1016/j.cam.2021.113912}}.}}

\begin{document}

\begin{abstract}
This paper introduces an efficient approach to solve quadratic and nonlinear programming problems subject to linear equality constraints via the Theory of Functional Connections.  This is done without using the traditional Lagrange multiplier technique.  More specifically, two distinct expressions (fully satisfying the equality constraints) are provided, to first solve the constrained quadratic programming problem as an unconstrained one for closed-form solution.  Such expressions are derived via using an optimization variable vector, which is called the free vector $\B{g}$ by the Theory of Functional Connections.  In the spirit of this Theory, for the equality constrained nonlinear programming problem, its solution is obtained by the Newton's method combining with elimination scheme in optimization.
Convergence analysis 
is supported by a numerical example for the proposed approach. 

\end{abstract}

\maketitle

\noindent \textbf{Keywords.} Quadratic and nonlinear programming; Linear equality constrained minimization problem; Newton's method and elimination scheme in optimization; Unconstrained minimization problem; Convergence; Theory of functional connections  

\vskip10pt

\noindent \textbf{Mathematics Subject Classification.} 65N99

\vfill

\noindent $^{*}$Corresponding author: \textit{Tina Mai}; $^a$Institute of Research and Development, Duy Tan University, Da Nang, 550000, Vietnam; $^b$Faculty of Natural Sciences, Duy Tan University, Da Nang, 550000, Vietnam; \textit{maitina@duytan.edu.vn}\\

\noindent $^c$\textit{Daniele Mortari}; Department of Aerospace Engineering, Texas A\&M University, College Station, Texas, USA; \textit{mortari@tamu.edu}


\newpage


\section{Introduction}\label{intro}


The \TFC\ (TFC)\footnote{This theory, initially called ``Theory of Connections,'' has been renamed for two reasons. First, the ``Theory of Connections'' already identifies a specific theory in differential geometry.  Second, what this theory is actually doing is ``Functional Interpolation'' as it provides \emph{all} functions satisfying a set of linear equality constraints in rectangular domains of the $n$-dimensional space.} is a recently proposed mathematical framework to perform \textit{functional interpolation}, which can be considered a generalization of interpolation. In particular, this methodology analytically derives functionals, called \textit{constrained expressions}, to represent \emph{all} functions satisfying a set of given $n$ linear equality constraints and to transform constrained problems into unconstrained ones, via the \textit{free function} $g (x)$ \cite{U-ToC}.  This function $g (x)$ was used in classical \textit{variable reduction method} (for generating null-space matrices \cite{2009proj}) as well as in classical \textit{null-space methods} (for quadratic programming (QP) \cite{QP3, 2009r}), and it has recently been named free function by the TFC.  Note that $g (x)$ can be discontinuous, partially defined, and even the Dirac delta function (as long as it is defined on where the constraints are defined). 

Consider a univariate function $y(x)$ subject to $n$ constraints, two formal and equivalent definitions of constrained expressions are introduced in \cite{U-ToC, M-ToC, PhD-Leake} as follows:
\begin{align}
	y (x) &= y(x, g(x)) = g (x) + \sum_{k = 1}^n \eta_k(x, g(x)) \ s_k (x)  \,, \label{eta} \\
	y (x) &= y(x, g(x)) = g (x) + \sum_{k = 1}^n \rho_k (x, g (x)) \ \phi_k (x, \B{s} (x))  \,. \label{phirho}
\end{align}
Here, $g(x)$ is a free function which can be arbitrarily chosen.  In \eqref{eta}, $\eta_k(x, g(x))$ are the \textit{functional coefficients} whose expressions are derived by imposing the given $n$ constraints on \eqref{eta}, 
and $\B{s} (x) = \{s_1 (x), \cdots, s_n (x)\}$ is a set of $n$ user-defined linearly independent \textit{support functions}.  In \eqref{phirho}, $\phi_k(x, \bfa{s}(x))$ are the \textit{switching functions} which imply changing between the constraints, and $\rho_k (x, g (x))$ are the \textit{projection functionals} which project the free function $g(x)$ to each provided $k$th constraint.  

For example, given constants $x_1, x_2, \beta_1, \beta_2\,,$ consider a function $y(x)$ subject to the constraints 
\begin{equation}\label{tie1}
\left.\dfrac{\dd y}{\dd x}\right|_{x = x_1} = \beta_1 \quad \text{and} \quad \left.\dfrac{\dd y}{\dd x}\right|_{x = x_2} = \beta_2\,.
\end{equation}
Using the form \eqref{eta} with $s_1(x) = x$ and $s_2(x) = x^2\,,$ the above constraints \eqref{tie1} can be expressed as
\begin{equation}\label{tie2}
y(x) = y(x, g(x)) = g(x) + \eta_1 \, s_1(x) + \eta_2 \, s_2(x) = g(x) + \eta_1 \, x + \eta_2 \, x^2\,,	
\end{equation}
where the two constants $\eta_1$ and $\eta_2$ are computed from \eqref{tie1}--\eqref{tie2} as follows:
\begin{align*}
&\begin{bmatrix} 
y'(x_1) - g'(x_1) \\ 
y'(x_2) - g'(x_2)	
\end{bmatrix} 
= 
\begin{bmatrix} 
	s_1'(x_1) & s_2'(x_1) \\ 
	s_1'(x_2) & s_2'(x_2)
\end{bmatrix} \,
\begin{bmatrix} 
	\eta_1 \\ 
	\eta_2
\end{bmatrix} 
\Leftrightarrow
\begin{bmatrix} 
	\eta_1 \\ 
	\eta_2
\end{bmatrix} 
=
\begin{bmatrix} 
	1 & 2x_1 \\ 
	1 & 2x_2
\end{bmatrix}^{-1} 
\, 
\begin{bmatrix} 
	\beta_1 - g'(x_1) \\ 
	\beta_2 - g'(x_2)
\end{bmatrix}
\\
\Leftrightarrow
&\begin{bmatrix} 
	\eta_1 \\ 
	\eta_2
\end{bmatrix} 
=
\frac{1}{2(x_2 - x_1)}
\begin{bmatrix} 
	2x_2 & -2x_1 \\ 
	-1   & 1
\end{bmatrix}
\, 
\begin{bmatrix} 
	\beta_1 - g'(x_1) \\ 
	\beta_2 - g'(x_2)
\end{bmatrix}
\\
\Leftrightarrow
&\begin{bmatrix} 
	\eta_1 \\ 
	\eta_2
\end{bmatrix} 
=
\begin{bmatrix} 
\dfrac{2x_2}{2(x_2-x_1)}(\beta_1 - g'(x_1)) + \dfrac{-2x_1}{2(x_2-x_1)}(\beta_2 - g'(x_2)) \\ 
\dfrac{-1}{2(x_2-x_1)}(\beta_1 - g'(x_1)) + \dfrac{1}{2(x_2-x_1)}(\beta_2 - g'(x_2))
\end{bmatrix}\,.
\end{align*}

Now, substituting the computed $\eta_1$ and $\eta_2$ into \eqref{tie2}, we derive the following expression of the form \eqref{phirho}:
\begin{equation}\label{tie2b}
	y(x) = y(x, g(x)) = g (x) +  \underbrace{(\beta_1 - g'(x_1))}_{\rho_1(x, g(x))} \,  \underbrace{\dfrac{x \, (2 x_2 - x)}{2 (x_2 - x_1)}}_{\phi_1 (x,\bfa{s}(x))} +
	 \underbrace{(\beta_2 - g'(x_2))}_{\rho_2(x, g(x))} \, 
	  \underbrace{\dfrac{x \, (x - 2 x_1)}{2 (x_2 - x_1)}}_{\phi_2 (x,\bfa{s}(x))}\,.
\end{equation}

This expression satisfies the constraints \eqref{tie1}, that is, $y'(x_1) = \beta_1$ and $y'(x_2) = \beta_2\,,$ \textit{no matter what the free function $g (x)$ is}.  Also, \eqref{tie2b} identifies the projection functionals $\rho_k (x, g (x))\,.$  Besides, the meaning of the switching functions $\phi_k (x) = \phi_k (x,\bfa{s}(x))$ in \eqref{tie2b} can be understood as follows: when the first constraint $y'(x_1) = \beta_1$ of \eqref{tie1} holds, we obtain from \eqref{tie2b} that $\phi'_1(x_1) = 1$ and $\phi'_2(x_1) = 0\,;$ when the second constraint $y'(x_2) = \beta_2$ holds, it follows that $\phi'_1(x_2) = 0$ and $\phi'_2(x_2) = 1\,.$	

\bigskip	

The Multivariate \TFC\ \cite{M-ToC} extends the original univariate theory \cite{U-ToC} to $n$ dimensions and to any-degree boundary constraints. This extension can be summarized by the expression
\begin{equation*}
	y (\B{x}) = h (c(\B{x})) + g (\B{x}) - h (g(\B{x}))\,,
\end{equation*}
where $\B{x} = (x_1, x_1, \dots, x_n)\T$ is the vector of $n$ orthogonal coordinates, $c (\B{x})$ is a function specifying the boundary constraints, $h (c(\B{x}))$ is \emph{any} interpolating function satisfying the boundary constraints, and $g (\B{x})$ is the free function. Several examples of constrained expressions can be found in \cite{U-ToC,M-ToC}. 

\bigskip

The TFC has been developed for point, derivative, integral, infinite, component constraints, and any linear combination of them for univariate functions \cite{U-ToC, PhD-Leake} as well as multivariate functions \cite{M-ToC, PhD-Leake}, in rectangular domains (with initial investigation) then in generic domain (as an extension) via domain mapping \cite{TFCbij}. The main feature of constrained expressions is that they allow restricting the search function space of a constrained problem to just the space of its feasible solutions
satisfying the constraints.  By this way, a large number of constrained optimization problems can be transformed into unconstrained ones, which can be solved by more simple, efficient, robust, reliable, fast, and accurate methods.

The first TFC application was in solving linear \cite{LDE} and nonlinear \cite{NDE} ODEs. This has been done by expanding the free function $g (x)$ in terms of a set of basis functions (for example, orthogonal polynomials, Fourier transforms, etc.). Linear or iterative nonlinear least-squares method is then used to find the coefficients of the expansion. This TFC approach for solving ODEs has many advantages over traditional methods: 1) it consists of a unified framework to solve IVP, BVP, or multi-valued problems, 2) it provides an analytically approximated solution that can be used for subsequent manipulation, 3) the solution is usually obtained in millisecond and at machine error accuracy, 4) the procedure is numerically robust (with small condition number), and 5) it can solve the ODEs subject to many different constraint types as above.  Additionally, this technique has recently been applied to solve a wide range of problems \cite{Selected}.

 Other specific TFC applications are as follows: homotopy continuation for control problems \cite{homotopyTFC}, epidemiological models \cite{TFCcovid}, radiative transfer problems \cite{TFCradi}, rarefied-gas dynamics \cite{TFCgas}, Timoshenko-Ehrenfest beam \cite{TFCfem}, hybrid systems \cite{TFChybrid}, machine learning \cite{TFCvectorML, TFCdeep, TFCextreme, TFCintercept}, orbit transfer and propagation \cite{JohnstonThesis, TFC2orbit, TFCporbit}, optimal control problems via indirect methods, relative motion \cite{TFCmotion}, landing on small and large planetary bodies \cite{TFCland}, and intercept problems \cite{TFCintercept}.	

The \textit{first purpose} of this paper is to complete the initial study presented in \cite{Selected}, which has shown that some classical optimization problems, such as constrained quadratic programming (QP), can be solved in closed-form and in efficient way by the \TFC (TFC).  For linear equality constrained QP, there are three classically basic approaches \cite{QP3, 2009r, fullsp08}: \textit{full-space, range-space, 
or null-space approaches}; and each of these can utilize either \textit{direct (factorization) method} or \textit{iterative (conjugate-gradient) method}.
In this paper, employing the TFC, two distinct approaches are introduced leading to unconstrained QP problems having closed-form solution, without using the classical Lagrange multiplier technique.  The \textit{first TFC approach} specifies the constrained expression \eqref{eta}.  This transforms the initial constrained optimization problem to an unconstrained one, where the variable is free vector $\B{g}\,.$  The resulting linear system is solved by using the Moore--Penrose pseudo-inverse.  The \textit{second TFC approach} takes advantage of the classical \textit{variable reduction method} to express the linear equality constraints.  In this case, the free vector $\B{g}$'s dimension is reduced.
Thus, benefiting the traditional \textit{null-space methods}, this second approach also transforms the constrained QP problem into an unconstrained linear system, but now it has nonsingular coefficient matrix.  Due to the reduction in sizes of both the free vector and the coefficient matrix, this second approach is much faster than the first one while giving the same exact solution.  To support these two approaches, numerical examples are provided.

For linear equality constrained nonlinear programming (NLP), there are several important classical approaches \cite{Boyd}.  The first approach (usually not possible) is \textit{analytically} finding a solution for its corresponding KKT system of equations.  The second approach is \textit{eliminating the equality constraints} to reduce the constrained problem to an unconstrained one (having fewer variables), which then can be handled by unconstrained optimization algorithms.  The third approach is solving the \textit{dual problem} utilizing an unconstrained optimization algorithm.  The fourth approach is the Newton's method with equality constraints (which is the most natural extension of the Newton's method without constraints): it finds solution for the NLP problem on the affine set of equality constraints.  The fifth method is a combination of the Newton's method (the fourth) and the elimination scheme (the second).  It is interesting that the iterates in Newton's method for the equality constrained NLP problem (the fourth) coincide with the iterates in Newton's method utilized to the unconstrained reduced problem (the fifth), and so does the convergence analysis.

Therefore, the \textit{second purpose} of this paper is to apply the TFC to the linear equality constrained convex NLP problem in the fifth method's spirit, that is, the Newton's method (by expanding the nonlinear objective function up to second order) combining with the elimination scheme toward an unconstrained minimization problem \cite{Boyd}.  Our main contribution is the convergence analysis supported by a numerical example.  In particular, we derive in the convergence analysis the termination criterion, quadratic rate,
and an upper bound on the total number of iterations required to attain a given accuracy.
For generality, we call this approach TFC because it can handle functional interpolation and it has a potential to extend to other algorithms (such as regular gradient descent methods, quasi-Newton and modified Newton methods) in addition to the classical Newton's method in this paper.   

It is crucial to note here that the TFC (eliminating constraints) is suitable for several classes of problems.  For example, the TFC is better than the finite element method (FEM) in some aspects \cite{TFCfem}.  Also, because the TFC finally leads to unconstrained NLP, gradient methods can be applied as well as accelerated by \textit{Picard-Mann} hybrid iterative process and so on \cite{psurvey}.  Mathematically, these advantages are  enjoyable.  In many cases \cite{Boyd}, however, it is better to keep the
equality constraints because eliminating them can make the problems more difficult to understand and analyze, or destroy the algorithm's efficiency. 
In those case, some method that directly works with the equality constraints (without eliminating them) can preserve the problem's structure.

\bigskip
 
The paper is organized as follows.  In Section \ref{2belts}, preliminaries and two TFC expressions of linearly under-determined equality constraints are presented with a numerical example for the second constrained expression's application.  Section \ref{qp0} is about introducing the quadratic programming problem subject to linear equality constraints and classical basic approaches to solve it.  We provide in Sections \ref{ibelt} and \ref{iibelt} the theory as well as numerical examples for solving quadratic programming with the first and second constrained expressions, respectively.  In Section \ref{nlps}, linear equality constrained nonlinear programming (NLP) is introduced together with the classical Newton's method.  We show in Section \ref{converge} the convergence analysis of the equality constrained NLP using the TFC (as a combination of the Newton's method and elimination scheme), with quadratic rate and bounded number of iterations.  In Section \ref{enlp}, a  numerical example of equality constrained NLP is provided to illustrate this convergence.  Section \ref{discuss} is regarding open discussion about accelerating the TFC by utilizing the Picard, mixed Picard-Newton, and Picard-Mann methods.  Conclusions are assembled in Section \ref{conclude}.


\section{Expressions of linearly under-determined equality constraints}\label{2belts}

First, we refer the readers to \cite{mcl} for the basic preliminaries.  The functions are denoted by italic lowercase letters (e.g., $f$), vector fields in $\mathbb{R}^d$ and matrix fields in $\mathbb{R}^{q \times s}$ are respectively represented by italic boldface lowercase and uppercase letters (e.g., $\bfa{b}$ and $\bfa{A}$).  
The sets of functions, vector fields in $\mathbb{R}^d\,,$ and matrix fields in $\mathbb{R}^{q \times s}$ are respectively denoted by italic, boldface, and special Roman capitals (e.g., $S\,,$ $\bfa{S}\,,$ and $\mathbb{S}$).  Throughout this paper, the symbol $\nabla$ stands for gradient.  The Euclidean norm of a vector $\bfa{g}$ is denoted by $\norm{ \bfa{g} }$ or $\norm{ \bfa{g} }_2\,.$  The Frobenius norm of a matrix $\B{A}$ is presented by $|\B{A}| : = \sqrt{\B{A} \cdot \B{A}} = \sqrt{\text{tr}(\B{A}\T \B{A})}\,.$

\vspace{10pt} 

Consider an under-determined linear system of equations
\begin{equation}\label{LinearSystem}
	\B{A} \, \B{x} = \B{b},
\end{equation}
where $\B{A} \in \mathbb{R}^{m\times n}$ ($m < n$) as well as $\B{b}\in \mathbb{R}^m$ are given, and $\B{x}\in \mathbb{R}^n$ is an unknown vector. We assume the consistent case that $\text{rank} (\B{A}) = m\,,$ so there is at least one
solution for \eqref{LinearSystem}.  In this paper, we do not consider the inconsistent case (with no solution) when the system \eqref{LinearSystem} forms distinctly parallel (not coinciding) hyperplanes. 

\subsection{First constrained expression}\label{1tie}

For the system \eqref{LinearSystem}, TFC is applied to specify the constrained expression form (\ref{eta}) as follows:
\begin{equation}\label{DM2}
	\B{x} (\B{g}) = \B{g} + \B{H} \, \B{\eta} (\B{g})\,,
\end{equation}
where $\B{g} \in \mathbb{R}^n$ is a free vector, $\B{H} \in \mathbb{R}^{n\times m}$ is an assigned support matrix with rank$(\B{H}) = m\,,$ 
and $\B{\eta}(\B{g}) \in \mathbb{R}^m$ is a coefficient vector that is a function of the free vector $\B{g}$. The expression of $\B{\eta}$ is derived by imposing the equality constraints \eqref{LinearSystem} on \eqref{DM2}. That is,
\begin{equation}\label{eta1}
	\B{A} (\B{g} + \B{H} \, \B{\eta}) = \B{b} \quad \Leftrightarrow \quad \B{A} \B{H} \B{\eta} = \B{b} - \B{A} \B{g} \quad \Leftrightarrow \quad \B{\eta} = (\B{A} \B{H})^{-1} \left(\B{b} - \B{A} \B{g}\right)\,.
\end{equation}
Substituting this expression of $\B{\eta}$ into  Eq.\ (\ref{DM2}), all solutions of the system (\ref{LinearSystem}) can be obtained in the form
\begin{equation}\label{1sta}
	\boxed{ \B{x} (\B{g}) = \B{g} + \B{H} \left(\B{A} \B{H}\right)^{-1} \left(\B{b} - \B{A} \, \B{g}\right)\,. }
\end{equation}
Equivalently, we get the following \textit{first constrained expression} for \eqref{LinearSystem}:
\begin{equation}\label{1st}
	\boxed{ \B{x} = \B{x}_0 + \B{D} \, \B{g}} \qquad {\rm where} \quad \begin{cases} \B{x}_0 = \B{H} (\B{A} \B{H})^{-1} \B{b} \\ \B{D} = \B{I}_{n\times n} - \B{H} (\B{A} \B{H})^{-1} \B{A} \,.\end{cases}
\end{equation}

The expression \eqref{1sta} or \eqref{1st} shows that $\B{H}$ can be \emph{any} matrix making the $m\times m$ matrix $\B{A H}$ nonsingular. This condition is satisfied if, for example, 
\begin{equation}\label{HUAT}
	\B{H} = \B{A}\T
\end{equation}
since  $\B{A}$ has full row rank $m\,.$  Note that with this choice \eqref{HUAT}, it is tempting to try to simplify \eqref{1sta} and \eqref{1st}; but since the given $\B{A}$ is not a square matrix, we can not have the inverse $\B{A}^{-1}$ nor  $(\B{A}\T \B{A})^{-1} = \B{A}^{-1} (\B{A}\T)^{-1}\,.$

Another option for $\B{H}$ (to minimize the computation time of matrix $\B{D}$ in \eqref{1st}) is to set
	\begin{equation}\label{Hmin}
		\B{H} = \begin{bmatrix} \B{I}_{m\times m} \\ 0_{(n-m)\times m}\end{bmatrix} \B{P}_{m\times m}\,,
	\end{equation}
	where $\B{P}_{m\times m}$ is chosen to be some $m\times m$ permutation matrix so that the $m \times m$ matrix $\B{AH}$ is formed by $m$ linearly independent columns of $\B{A}$ dictated by $\B{P}$ \cite{2009proj}.
\subsection{Second constrained expression}\label{2tie}

Now, we consider the case (\ref{HUAT}).  If $\B{g} = \B{0}$, then using the Moore--Penrose right inverse
\begin{equation}\label{mppr}
	\B{A}^{+} = \B{A}\T \, \left(\B{A} \, \B{A}\T\right)^{-1} \,, 
\end{equation}	
Eq.\ (\ref{1sta}) has the following solution:
\begin{equation}\label{Moore-Penrose}
	\B{x}_0 = \B{x} (\B{0}) =  \B{A}\T \, \left(\B{A} \, \B{A}\T\right)^{-1} \, \B{b} = \B{A}^{+} \B{b}\,.
\end{equation}
It thus follows from \eqref{1sta} that the difference $\Delta \B{x} = \left[\B{I} - \B{A}\T \, \left(\B{A} \, \B{A}\T\right)^{-1} \B{A}\right] \, \B{g}$ identifies the term allowing us to obtain all solutions which are different from the solution \eqref{Moore-Penrose}. 
In this spirit, Eq.\ (\ref{1sta}) can be written as
\begin{equation}\label{DCE}
	\B{x} (\B{g}) = \B{x}_0 + \B{D} \, \B{g} \qquad \text{where} \qquad \B{D} = \B{I}_{n \times n} - \B{A}\T \, \left(\B{A} \, \B{A}\T\right)^{-1} \B{A}\,.
\end{equation}
Premultiplying this $\B{D}$ by $\B{A}$, we obtain
\begin{equation*}
	\B{A} \, \B{D} = \B{A} - \B{A} \, \B{A}\T \, \left(\B{A} \, \B{A}\T\right)^{-1} \B{A} = \B{A} - \B{A} = \B{0}\,.
\end{equation*}
In \eqref{DCE}, the $n \times n$ matrix $\B{D} = \B{I}_{n \times n} - \B{A}\T \, \left(\B{A} \, \B{A}\T\right)^{-1} \B{A}$ is called the orthogonal projection matrix onto the nullspace of $\B{A}\,.$ Thus, $\B{D}^2 = \B{D}$ and $\B{D}$ is also a null-space matrix for $\B{A}$ (that is, $\B{D}\B{g} \in \tu{null}(\B{A})$ for any free vector $\B{g} \in \mathbb{R}^{n}$) \cite{2009proj}.
	
Let $\B{N} \in \mathbb{R}^{n \times (n-m)}$ be a matrix having columns as a \textit{basis} for the nullspace of matrix $\B{A} \in \mathbb{R}^{m \times n}$ \cite{2009proj}.  For brevity, a basis for the null space of $\B{A}$ is call a \textit{null-space basis} \cite{QP3} or \textit{null basis} \cite{null2alg}, and a column of a null basis is called a \text{null vector} for $\B{A}$ \cite{null2alg}.
Automatically, $\B{A N} = \B{0}_{m\times (n-m)}\,.$  If this \textit{null basis} $\B{N}$ replaces $\B{D}$ in \eqref{DCE}, then the free vector $\B{g}$'s size (instead of $n$) can be reduced to the minimum ($n-m$), and thus any feasible vector $\B{x}(\B{g})$ has the form 
\begin{equation}\label{2nd}
	\B{x} (\B{g}) = \B{x}_0 + \B{N} \B{g}\,,
\end{equation}
which is call the TFC \textit{second constrained expression} for \eqref{LinearSystem}.  Here, $\B{x}_0$ is \textit{any} solution vector of the system (\ref{LinearSystem}), and $\B{g} \in \mathbb{R}^{n-m}$ is the free vector.  Note that $\B{N}\B{g} \in \tu{null}(\B{A})$ and $\B{A} (\B{N}\B{g}) = \bfa{0}_{m\times 1}\,.$  The expression \eqref{2nd} was introduced in classical \textit{variable reduction method} in generating null-space matrices \cite{2009proj} and also was used in traditional \textit{null-space methods} for QP \cite{QP3, 2009r}.

The proof that expression (\ref{2nd}) (or similar (\ref{DCE})) provides \textit{all} solutions of the system (\ref{LinearSystem}) is immediate.  Indeed, let $\bfa{g} \in \mathbb{R}^{n-m}$ be the free vector.  Then, $\B{N} \bfa{g} \in \text{Null}(\B{A})\,.$ Thus, $\B{A}(\bfa{x}_0 + \B{N} \bfa{g}) = \B{A}(\bfa{x}_0) + \B{A}(\B{N} \bfa{g}) = \B{A}(\bfa{x}_0) + \bfa{0} = \bfa{b}\,.$  Conversely, let $\B{A}\bfa{x} = \bfa{b}\,.$  Then, $\bfa{x} = \bfa{x}_0 + \bfa{x} - \bfa{x}_0\,.$  We have $\bfa{b} = \B{A}\bfa{x} = \B{A}(\bfa{x}_0 + \bfa{x} - \bfa{x}_0) = \B{A} \bfa{x}_0 + \B{A}(\bfa{x} - \bfa{x}_0)\,.$  Since $\B{A}\bfa{x}_0 = \bfa{b}\,,$ it follows that $\B{A}(\bfa{x} - \bfa{x}_0) = \bfa{0}\,,$ implying $\bfa{x} - \bfa{x}_0 = \B{N}\bfa{g} \in \text{Null}(\B{A})$ for any $\bfa{g} \in \mathbb{R}^{n-m}\,.$


There is a parallelism between the interpolation and problem of solving an under-determined linear system of equations. For example, as in (\ref{Moore-Penrose}), one solution $\B{x}_0$ (with the smallest $l^2$-norm value from all admissible $\B{x}$) 
for the under-determined linear system of equations (\ref{LinearSystem}) is given by the Moore--Penrose right inverse \eqref{mppr}.  That is \cite{mpp17},
\begin{equation*}
||\B{x}_0||_2 = || \B{A}^{+} \B{A} \bfa{x}||_2 \leq ||\B{u}||_2 \,, \quad \tu{for all } \B{u} \tu{ satisfying } \B{A}\bfa{u} = \B{A}\bfa{x} = \bfa{b}\,.	
\end{equation*}	 
%

\noindent Many other feasible solutions $\B{x}_0$ can be obtained by assigning random values to the $(n - m)$ variables and then solving \eqref{LinearSystem} for the remaining $m$ unknowns. Using \textit{any} existing solution $\B{x}_0$ of (\ref{LinearSystem}), the TFC can then provide \textit{all} possible solutions. An equivalent approach was done for interpolation is that using any existing interpolation function, the TFC can then provide all interpolating functions. 

\subsection{Discussion on the two constrained expressions}\label{I-IIbelts}

There is an important difference between the use of (\ref{DCE}) and (\ref{2nd}) expressions. At first glance, the null-space formulation (\ref{2nd}) appears to be more convenient because the free vector has smaller dimensions ($n-m < n$). However, \cite{nullcomplex} 
has shown that the null-space problem (finding a sparsest null-space basis matrix $\B{N}$ for $\B{A}$) is NP-hard because finding a sparsest null vector of $\B{A}$ is NP-complete. On the other hand, by setting the support matrix $\B{H}$ as $\B{H} = \B{A}\T$ in \eqref{HUAT}, the execution of (\ref{DCE}) has complexity in polynomial time via naive matrix multiplication and a new Gauss--Jordan elimination method \cite{mpptime} for the Moore--Penrose
right inverse \eqref{mppr}.
%
 Choosing $\B{H}$ as in (\ref{Hmin}), the resulting complexity for performing (\ref{DCE}) is even smaller than using \eqref{HUAT}.  Since (\ref{DCE}) can be executed in polynomial time which is at most $\mathcal{O} \left(n^3\right)$, it follows that (\ref{DCE}) has lower complexity than \eqref{2nd}.


\subsection{An application of the second constrained expression}\label{zerofind} 
In the spirit of second equality constrained expression \eqref{2nd}, the following algebraic zero-finding problem was solved by TFC homotopy continuation in \cite{homotopyTFC}:
\begin{equation}\label{HomotopyExample}
\bfa{f}(\B{x}) = \B{0}\,,
\end{equation}
where $\B{x} = (x_1,x_2) \in \mathbb{R}^2$ and $\B{f} \colon \mathbb{R}^2 \to \mathbb{R}^2$ is the $\mathcal{C}^2$ 
vector-valued function
\begin{equation*}
	\B{f} (x_1, x_2) = \begin{bmatrix} a (x_1 + x_2) \\
		a (x_1 + x_2) + (x_1 - x_2) \left((x_1 - c)^2 + x_2^2 - d\right)\end{bmatrix}\,,
\end{equation*}
with $a = 4$, $c = 2$, and $d = 1\,.$

All solutions of this problem  \eqref{HomotopyExample} can be found by using (\ref{2nd}) because the first equation of \eqref{HomotopyExample} can be seen as
\begin{equation}\label{0eg}
	\B{A} \, \B{x} =
	\begin{bmatrix} 
			a & a\\
			0 & 0	
	\end{bmatrix}  \begin{bmatrix} x_1 \\ x_2\end{bmatrix} = \begin{bmatrix} 0 \\ 0\end{bmatrix} = \B{b} \,.
\end{equation}
An orthonormal basis for the nullspace of $\B{A}$ is 
\[\B{N} = \dfrac{\sqrt{2}}{2} \begin{bmatrix} -1 \\ +1\end{bmatrix}\,.\]

Since $\B{b} = \B{0}$, it follows from \eqref{Moore-Penrose} that $\B{x}_0 = \B{0}\,.$ Thus, all solutions to \eqref{0eg} are of the form $\B{x} = \B{N} g$, where $g$ is a scalar variable.  Substituting this form of $\B{x}$ into the second equation of \eqref{HomotopyExample}, we obtain and solve an algebraic cubic equation in the variable $g$ to get all the three solutions of \eqref{HomotopyExample} as follows:
\begin{equation*}
	\B{x}_1 = \begin{bmatrix} 0 \\ 0\end{bmatrix}\,, 
	\qquad
	\B{x}_{2} = \left(1 + i \dfrac{\sqrt{2}}{2}\right) \begin{bmatrix} 1 \\ -1\end{bmatrix}\,,
	\quad 
	\B{x}_{3} = \left(1 - i \dfrac{\sqrt{2}}{2}\right) \begin{bmatrix} 1 \\ -1\end{bmatrix}\,.
\end{equation*}
	
We wish to emphasize in this example that if rank$(\B{A}) = n - 1$ (where $n$ is the number of variables), then the system \eqref{HomotopyExample} becomes univariate and can be solved by zero-finding methods.

\section{Quadratic programming subject to linear equality constraints}\label{qp0}

The problem is to find the extreme of quadratic function $f (\B{x})$, subject to $m < n$ linear constraints:
\begin{equation}\label{DM1}
   {\rm\bf QP:} \; \minx_{\B{x} \in \mathbb{R}^n} \; \left\{f (\B{x}) = \dfrac{1}{2} \, \B{x}\T \B{Q} \B{x} + \B{c}\T \B{x} \right\}\,,   \quad \B{A} \, \B{x} = \B{b}\,,
\end{equation}
where $\bfa{x} \in \mathbb{R}^n$ is an unknown vector, and we are given $f (\B{x}): \mathbb{R}^n \to \mathbb{R}$, $\B{Q} = \B{Q}\T \in \mathbb{R}^{n\times n}$, $\B{c} \in \mathbb{R}^n$, $\B{A} \in \mathbb{R}^{m\times n}$, $\B{b} \in \mathbb{R}^m$, and rank$(\B{A}) = p \le m$.  Note that $\B{Q}$ can always be presumed to be symmetric.  Indeed, 
	\[\B{x}\T \B{Q} \B{x} = \frac{1}{2} \B{x}\T (\B{Q} + \B{Q}\T) \B{x}\,,\]
	 so $\B{Q}$ could be replaced by the symmetric matrix \[\bar{\B{Q}} = \frac{1}{2}   (\B{Q} + \B{Q}\T) \,.\]
	 


Most of the current approaches to solve the problem (\ref{DM1}) can be found in \cite{QP3, 2009r, fullsp08} supported by \cite{QP1,QP2,QP4,QP5}. The classical approach uses the Lagrange multiplier technique, which leads to an equivalent linear system to \eqref{DM1}, namely the Karush-Kuhn-Tucker (KKT) system:
\begin{equation}\label{KKT}
    \begin{bmatrix} \B{Q} & - \B{A}\T \\ \B{A} & \B{0}_{m\times m}\end{bmatrix} \begin{Bmatrix} \B{x} \\ \B{\lambda}\end{Bmatrix} = \begin{Bmatrix*}[r] -\B{c} \\ \B{b}\end{Bmatrix*}\,,
\end{equation}
where $\B{\lambda}$ is a vector containing $m$ Lagrange multipliers. This linear system \eqref{KKT} provides a unique solution for \eqref{DM1} if the following conditions hold \cite{QP3}: 
\noindent 1) The matrix $\B{A}$ has full rank $m\,,$ and 
\noindent 2) the reduced Hessian matrix $\B{N}\T \B{Q} \B{N}$ is positive-definite, where $\B{N} \in \mathbb{R}^{n \times (n-m)}$ is a matrix whose columns form a basis for the nullspace of $\B{A}$ (that is $\B{N}$ has full rank and $\B{AN} = \bfa{0}_{m \times (n-m)}$).

\bigskip

The problem \eqref{DM1} can be solved by three classically basic approaches \cite{QP3, 2009r, fullsp08}: \textit{full-space, range-space, 
and null-space approaches}; each of these can use either \textit{direct (factorization) method} or \textit{iterative (conjugate-gradient) method}.  

More specifically, full-space methods \cite{QP3, QP3, 2009r} directly handle the indefinite KKT system \eqref{KKT} by employing, for example, a triangular factorization which is usually in the form of an efficient \textit{symmetric indefinite factorization}.  


Range-space methods (the so-called Schur-complement method or block-elimination approach) \cite{2009r, QP3, fullsp08} require $Q$ to be positive-definite (and thus nonsingular) so that $\bfa{x}$ can be obtained by eliminating the first block of \eqref{KKT}:
\begin{align*}
	\B{A} \B{Q}^{-1}\B{A}\T \B{\lambda} &= \B{A} \B{Q}^{-1} \B{c} + \B{b}\,,\\
	\B{Q}\B{x} &= \B{A}\T \B{\lambda} - \B{c}\,.
\end{align*}

Null-space methods (which form the second TFC constrained expression to be discussed below) use an $n \times (n-m)$ matrix $\B{N}$ which is a basis for the nullspace of $\B{A}$ \cite{fullsp08}.  These methods do not require that $Q$ is nonsingular and hence has wider applicability than the range-space methods \cite{QP3}.

\bigskip

In the following Sections \ref{ibelt} and \ref{iibelt},
 two distinct approaches are introduced using the two TFC constrained expressions presented in Subsections \ref{1tie} and \ref{2tie}. 
  The \textit{first approach} employs the TFC to specify constrained expression \eqref{eta} whose form is in Subsection \ref{1tie}, which uses a free vector $\B{g}$ having the same size as the solution vector $\B{x}$'s dimension. This approach transforms the given linear equality constrained QP problem to an unconstrained linear system with (either nonsingular or singular) $n\times n$ coefficient matrix. The solution of this system is obtained through the Moore--Penrose pseudo-inverse. The second approach utilizes the \textit{null-space methods} \cite{QP3, 2009r}, whose constrained expression is presented in Subsection \ref{2tie} using \textit{variable reduction method} \cite{2009proj}.  This approach takes advantage of the basis $\B{N}$ for the nullspace of matrix $\B{A}$ and uses a \textit{free function} $\B{g}$ with dimension $(n-p)$, where $p = \text{rank}(\B{A})$.  Here, the constrained QP problem is also transformed into a linear system, but now it has nonsingular coefficient matrix.  Both approaches assume rank$(\B{A}) = p \le m$, where $m$ is the number of equality constraints.
  

\section{Quadratic programming with the first constrained expression}\label{ibelt}

For the constraint $\B{A} \, \B{x} = \B{b}$ associated with the problem (\ref{DM1}), the TFC is applied to specify constrained expression \eqref{eta} as in Section \ref{2belts}, that is \eqref{DM2}\footnote{This approach appeared in \cite{Selected} as one of the new TFC applications, which are different from solving differential equations.}:
\begin{equation*}
    	\B{x} (\B{g}) = \B{g} + \B{H} \, \B{\eta} (\B{g})\,.
\end{equation*}
Here, we are given $\B{A} \in \mathbb{R}^{m\times n}$ ($m < n$) with $\tu{rank}(\B{A}) = m$ and  $\B{b}\in \mathbb{R}^m\,,$  while $\B{x}\in \mathbb{R}^n$ is an unknown vector; $\B{H} \in \mathbb{R}^{n\times m}$ is an assigned matrix with rank$(\B{H}) = m$, $\B{g} \in \mathbb{R}^n$ is a free vector, and $\B{\eta} \in \mathbb{R}^m$ is a vector-valued function of $\B{g}\,.$ 

In Section \ref{2belts}, after computing the above $\B{\eta}\,,$ an equivalent form of \eqref{DM2} is \eqref{1st}:
\begin{equation*}
	\boxed{ \B{x} = \B{x}_0 + \B{D} \, \B{g}} \qquad {\rm where} \quad \begin{cases} \B{x}_0 = \B{H} (\B{A} \B{H})^{-1} \B{b} \\ \B{D} = \B{I}_{n\times n} - \B{H} (\B{A H})^{-1} \B{A}\,.\end{cases}
\end{equation*}
It follows from this expression that $\B{H}$ can be \emph{any} matrix, such as \eqref{HUAT} or \eqref{Hmin}, making the $m\times m$ matrix $\B{A H}$ nonsingular.

Now, substituting the above expression of $\B{x}$ from (\ref{1st}) into the problem (\ref{DM1}), we obtain the objective function $f(\B{x})$ in the form
\begin{equation*}
    h (\B{g}) = \dfrac{1}{2} (\B{x}_0 + \B{D} \, \B{g})\T \B{Q} (\B{x}_0 + \B{D} \, \B{g}) + \B{c}\T (\B{x}_0 + \B{D} \, \B{g}).
\end{equation*}
The problem (\ref{DM1}) thus becomes unconstrained. Its first-order optimal condition
\begin{equation*}
	\nabla h(\B{g}) = \B{0}
\end{equation*}
implies
\begin{equation}\label{DM4}
    \boxed{ \B{\cal A} \, \B{g} + \B{d} = \B{0}_{n\times 1}} \qquad \text{where} \quad \begin{cases} \B{\cal A} = \B{D}\T \B{Q} \B{D}\,, \\ \B{d} = \B{D}\T \B{Q} \B{x}_0 + \B{D}\T \B{c}\,.
    \end{cases}
\end{equation}
The resulting matrix $\B{\cal A}$ can be either nonsingular or singular because of $\B{D}$.  Indeed, it can happen that $\tu{det}(\B{D}) = \tu{det}(\B{I} - \B{H} (\B{A H})^{-1} \B{A}) \neq \B{0}$.  Hence, $\B{D}$ can be either nonsingular or singular (special case).  Also, we note that $\B{A}$ only has right inverses $\B{A}_R^{-1}$ with $\B{A A}_R^{-1} = \B{I}\,,$ and the nicest one of these $\B{A}_R^{-1}$ is the unique Moore--Penrose right inverse $\B{A}^{+} = \B{A}\T \, \left(\B{A} \, \B{A}\T\right)^{-1}\,.$ 
%
Similarly, $\B{H}$ only has left inverses $\B{H}_L^{-1}$ such that $\B{H}_L^{-1} \B{H} = \B{I}\,.$  Therefore, $\B{\cal A}$ can be either nonsingular or singular (special case).  Anyhow, regardless of the rank of $\B{\cal A}$ and regardless whether $\B{\cal A}$ is nonsingular or singular, all solutions (if any) of the system \eqref{DM4} are obtained by using the Moore--Penrose (generalized) pseudo-inverse $\B{\cal A}^{+}$ of $\B{\cal A}$ as follows:
\[ \B{g} =-\B{\cal A}^{+} \B{d} + (\B{I}_{n \times n} - \B{\cal A}^{+} \B{\cal A})\B{w}\,,\]
for arbitrary vector $\B{w} \in \mathbb{R}^n\,.$  In general (even when the system \eqref{DM4} is inconstent), $\B{g} = -\B{\mathcal A}^{+} \B{d}$ represents the least square (approximate) solution with minimum Euclidean norm $\| \cdot \|_2$ among all (approximate) solutions to  the system \eqref{DM4}. 
A computationally easy and exact method to compute the pseudo-inverse $\B{\mathcal A}^{+}$ of $\B{\cal A}$ is via applying the singular value decomposition (SVD).  Now, let $\B{\cal A}=\B{U\Sigma V}\Ts$ be the singular value decomposition of $\B{\cal A}$, then $\B{\cal A}^{+} = \B{V\Sigma^{+}U}\Ts$.  For such a square diagonal matrix $\B{\Sigma}$, its pseudo-inverse $\B{\Sigma}^{+}$ is obtained by taking the inverse of every nonzero element and keeping the zeros in position.  Note that $\B{U},\B{V}$ are orthogonal matrices.  Using \eqref{1st} and \eqref{DM4}, the solution of problem (\ref{DM1}) is then
\begin{equation}\label{DM}
	\boxed{ \B{x} 
		= \B{x}_0 - \B{D} \B{\mathcal A}^{+} \B{d} = \B{x}_0 - \B{D}\left(\B{D}\T \B{Q} \B{D}\right)^{+} \B{D}\T\left(\B{Q} \, \B{x}_0 + \B{c}\right)  =  \B{x}_0 - \B{D} \, \B{V\Sigma^{+}U}\Ts \, \B{d} } \,.
\end{equation}
Note that the analytical solution provided by (\ref{DM}) does not require $\B{Q}$ (from \eqref{DM1}) nor the Hessian $\B{N}\T \B{Q} \B{N}$ to be positive-definite, where $\B{N}$ is a null-space basis matrix for $\B{A}\,.$ 

\subsection{Equivalent reduced equality constraints}\label{ibeltr}

If rank$(\B{A}) = p < m$, then the matrix $\B{A}$ has only $p$ linearly independent rows. 
In addition, if there are no distinctly parallel hyperplanes from $\B{A} \, \B{x} = \B{b}$ (so no inconsistency), 
then
the constraint $\B{A} \, \B{x} = \B{b}$ can be transformed into an equivalent $\tilde{\B{A}} \, \B{x} = \tilde{\B{b}}$ system, where $\tilde{\B{b}} \in \mathbb{R}^p$ and $\tilde{\B{A}} \in \mathbb{R}^{p\times n}$ with rank$(\tilde{\B{A}}) = p$.  This transformation can be obtained via the rank revealing QR (RRQR) decomposition (\cite{rrqr3, rrqr4, rrqr1, rrqrn4}), which computes a decomposition of a matrix $\B{A} \in \mathbb{R}^{m \times n} \ (m < n)$ of the form (\cite{rrqr}):
\begin{equation}\label{rrqr}
 \B{AP} = \B{QR} = \B{Q} 
 \begin{bmatrix}
  \B{R}_{11} & \B{R}_{12}\\
  \B{0}_{(m-k)\times k} & \B{R}_{22}
 \end{bmatrix}\, \qquad \Rightarrow
  \qquad  \B{Q}\T \B{A} \, \B{x} = \B{RP}\T \, \B{x} = \B{Q}\T \, \B{b} \,,
\end{equation}
where $\B{Q}\in SO(m)$, $\B{R} \in \mathbb{R}^{m \times n}$ is upper trapezoidal such that $\B{R}_{11} \in \mathbb{R}^{k\times k}$ is upper triangular, $\B{R}_{12} \in \mathbb{R}^{k\times (n-k)}$, and $\B{R}_{22} \in \mathbb{R}^{(m-k)\times (n-k)}$.  The column permutation matrix $\B{P} \in \mathbb{R}^{n \times n}$ and the integer $k$ are selected such that $\| \B{R}_{22} \|_2$ is small (that is, $\B{R}_{22}$ is considered to be $\B{0}$) and $\B{R}_{11}$ is well-conditioned as well as possessing non-decreasing diagonal elements.  This decomposition was
introduced in \cite{rrqr1}, and the first algorithm to calculate it was suggested in \cite{rrqr2} as well as
thanks to the QR decomposition with column pivoting \cite{rrqrn1,rrqrn2, rrqrn3, rrqrn4} (at lower computational cost than a singular value decomposition).  

Computationally, setting $\B{t} =$ diag$(\B{R})$, where $|\B{t}_1| \ge |\B{t}_2|, \cdots$, we get the rank $p =k$ of $\B{A}$ as the maximum integer $k$ satisfying
\begin{equation}\label{kmax}
    |\B{t}_k| > \varepsilon \, \max \{m,n \} \, |\B{t}_1|\,,
\end{equation}
where $\varepsilon$ is a very small tolerance. 


Therefore, the equivalent reduced system is obtained by selecting the first $p = k$ rows of the system 
$\B{R P}\T \, \B{x} = \B{Q}\T \, \B{b}$, that is,
\begin{equation*}
    \tilde{\B{A}} \, \B{x} = \tilde{\B{b}}.
\end{equation*}
Hence, (\ref{1st}) becomes
\begin{equation*}
    \boxed{ \B{x} = \B{x}_0 + \B{D} \, \B{g}} \qquad {\rm where} \quad \begin{cases} \B{x}_0 = \tilde{\B{H}} (\tilde{\B{A}} \tilde{\B{H}})^{-1} \B{b} \\ \B{D} = \B{I}_{n\times n} - \tilde{\B{H}} (\tilde{\B{A}} \tilde{\B{H}})^{-1} \tilde{\B{A}}\,,
    \end{cases}
\end{equation*}
where $\tilde{\B{H}}$ can be set as $\tilde{\B{H}} = \tilde{\B{A}}\T$, and the solution follows from the previously derived approach for the case rank$(\B{A}) = m$.


\subsection{Monte Carlo tests for QP with first constrained expression}\label{nIbelt} 

For the problem \eqref{DM1}, we employ the formulas $\B{H} = \B{A}\T$ from \eqref{HUAT}, $\B{D} = \B{I}_{n\times n} - \B{A}\T \, \left(\B{A} \, \B{A}\T\right)^{-1} \B{A}$ from \eqref{DCE}, and $\B{x}_0 =  \B{A}\T \, \left(\B{A} \, \B{A}\T\right)^{-1} \, \B{b} = \B{A}^{+} \B{b}$ from (\ref{Moore-Penrose}).  Recall that the TFC solution to the problem \eqref{DM1} is \eqref{DM}:
\begin{equation*}
	\B{x} = \B{x}_0 - \B{D}\left(\B{D}\T \B{Q} \B{D}\right)^{+} \B{D}\T\left(\B{Q} \, \B{x}_0 + \B{c}\right)\,.
\end{equation*}

Using arrays $\B{A}$, $\B{b}$, $\B{Q}$, $\B{c}$ (from \eqref{DM1}) containing real floating-point numbers in $[-1, 1]$ drawn from a uniform distribution, we carry out a series of $10000$ Monte Carlo tests to compare the results provided by the TFC solution formula \eqref{DM} with the results obtained by MATLAB's \texttt{quadprog} built-in function for the problem \eqref{DM1}.  Note that in MATLAB's \texttt{quadprog}, we can choose either interior-point-convex or  trust-region-reflective algorithms, which are classical methods to solve \eqref{DM1} as well.  Therefore, these Monte Carlo tests also some how compare the TFC with those classical methods built in MATLAB's \texttt{quadprog}.  The results are shown in Table \ref{speed1}, having the array dimensions $n$ and $m$, the mean of elapsed time $\Delta t$ (in msec) and the mean of accuracy $e_2 = || \B{A} \, \B{x} - \B{b}||_2$ (which demonstrate that the computed TFC solution \eqref{DM} satisfies the given equality constraint in \eqref{DM1}).

\begin{table}[ht]
	\centering
	\begin{tabular}{|c|c|c|c|c|c|} \hline
		$n$ & $m$ & $\Delta t$ (TFC) & $\Delta t$ (MATLAB) & mean $e_2$ (TFC) & mean $e_2$ (MATLAB) \\ \hline \hline
		$10$ & $2$ & $0.21344$ & $3.2137$ & $1.3397\times 10^{-13}$ & $4.8483$ \\
		$10$ & $4$ & $0.18705$ & $3.2706$ & $6.4052\times 10^{-13}$ & $7.1054$ \\
		$10$ & $8$ & $0.18163$ & $3.3299$ & $9.328\times 10^{-10}$ & $8.9608$ \\
		$20$ & $4$ & $0.1916$ & $3.3255$ & $6.9779\times 10^{-13}$ & $9.8271$ \\
		$20$ & $8$ & $0.23446$ & $3.7902$ & $1.0836\times 10^{-12}$ & $14.4824$ \\
		$20$ & $16$ & $0.2619$ & $5.1378$ & $1.1712\times 10^{-10}$ & $20.6576$ \\
		$40$ & $8$ & $0.31959$ & $5.2355$ & $1.9041\times 10^{-12}$ & $20.366$ \\
		$40$ & $16$ & $0.29341$ & $3.7751$ & $8.2595\times 10^{-11}$ & $29.149$ \\
		$40$ & $32$ & $0.32228$ & $3.7274$ & $3.7067\times 10^{-10}$ & $41.4705$ \\
		$80$ & $16$ & $0.84889$ & $4.2288$ & $1.6493\times 10^{-11}$ & $40.814$ \\
		$80$ & $32$ & $0.94851$ & $5.1808$ & $3.4551\times 10^{-11}$ & $58.28$ \\
		$80$ & $64$ & $1.011$ & $5.8075$ & $5.7209\times 10^{-10}$ & $83.1533$ \\        \hline
	\end{tabular}
	\caption{Speed tests for QP with the first constrained expression.}
	\label{speed1}
\end{table}


\section{Quadratic programming with the second constrained expression}\label{iibelt}



We recall from Section \ref{2belts} that for the constraint $\B{A} \, \B{x} = \B{b}$ provided in the problem (\ref{DM1}), TFC derives a second constrained expression \eqref{2nd} (which was introduced in classical variable reduction method \cite{2009proj} and was used in traditional null-space methods for QP \cite{QP3, 2009r}):
\begin{equation}\label{pippo}
    \B{x} = \B{x}_0 + \B{N} \B{g} \qquad \text{with} \qquad \B{x}_0 = \B{A}\T (\B{A A}\T)^{-1} \,,
\end{equation}
where $\B{N} \in\mathbb{R}^{n\times (n-m)}$ is a matrix whose columns form a basis for the nullspace of matrix $\B{A} \in \mathbb{R}^{m\times n}$ (that is, $\B{A N} = \B{0}_{m\times (n - m)}$), $\tu{rank}(\B{A}) = m\,,$ and $\B{g} \in \mathbb{R}^{n-m}$ is the free vector. Such a null basis $\B{N}$ for $\B{A}$ can be computed in various ways: by SVD, QR, Cholesky decomposition as well as Gaussian elimination and Reduced Row Echelon Form.  Note that the dimension of the free vector $\B{g}$ in \eqref{pippo} is $(n-m)$, which is smaller than the size $n$ introduced in (\ref{DM2}) or \eqref{1st}.

Substituting this expression of $\B{x}$ for the one in problem (\ref{DM1}), we obtain $f(\B{x})$ in the form
\begin{equation*}
    h (\B{g}) = \dfrac{1}{2} (\B{x}_0 + \B{N} \, \B{g})\T \B{Q} (\B{x}_0 + \B{N} \, \B{g}) + \B{c}\T (\B{x}_0 + \B{N} \, \B{g})\,,
\end{equation*}
that is,
\begin{equation*}
    h (\B{g}) = \frac{1}{2}\B{g}\T \B{B} \, \B{g} + \left(\frac{1}{2} \B{x}_0\T \B{Q N} + \B{c}\T \B{N} \right) \B{g} + \frac{1}{2}\B{g}\T \B{N}\T \B{Q} \B{x}_0 + \frac{1}{2}\B{x}_0\T \B{Q} \B{x}_0 + \B{c}\T \B{x}_0\,.
\end{equation*}
Let
	\begin{equation*}
		\B{B} = \B{N}\Ts \B{Q N} \qquad \text{and} \qquad \B{p} = \B{N}\Ts \B{Q} \B{x}_0 + \B{N}\Ts \B{c}\,.
	\end{equation*}
	Note that the following gradients are with respect to $\B{g}\,:$
	\begin{equation*}
		\nabla (\B{g}\T \B{B} \, \B{g}) = (\B{B} + \B{B}\T) \B{g}\,, \qquad 
		\nabla(\B{y}\T \B{g}) = \B{y}\,, \qquad \text{and} \qquad 
		\nabla(\B{g}\T \B{u}) = \B{u}\,.
	\end{equation*}
Then the optimal problem is similar to the one in (\ref{DM}):
\begin{equation}\label{nulleq}
	\nabla h(\B{g}) = \B{0} \quad \Leftrightarrow \quad \B{B} \, \B{g} + \B{p} = \B{0}_{(n-m)\times 1}\,,
\end{equation}
but the only difference is that we now assume that the reduced Hessian $\B{B}$ is positive-definite, so $\B{B}$ is an $(n-m) \times (n-m)$ nonsingular matrix. 
Thus, the solution of \eqref{nulleq} is as follows \cite{null02e}:
\begin{equation*}
    \boxed{ \B{x} = \B{x}_0 - \B{N} \, \B{B}^{-1} \, \B{p} }\,,
\end{equation*}
which does not require nonsingularity of $\B{Q}$ \cite{QP3}.

If rank$(\B{A}) = p < m$, then the equality constraint $\B{A} \, \B{x} = \B{b}$ is transformed into the equivalent reduced equality constraint $\tilde{\B{A}} \, \B{x} = \tilde{\B{b}}$ as in Subsection \ref{ibeltr}.


\subsection{Numerical validation tests for QP with second constrained expression}\label{2ndne}

\begin{figure}[h]
    \centering\includegraphics[width=1\linewidth]{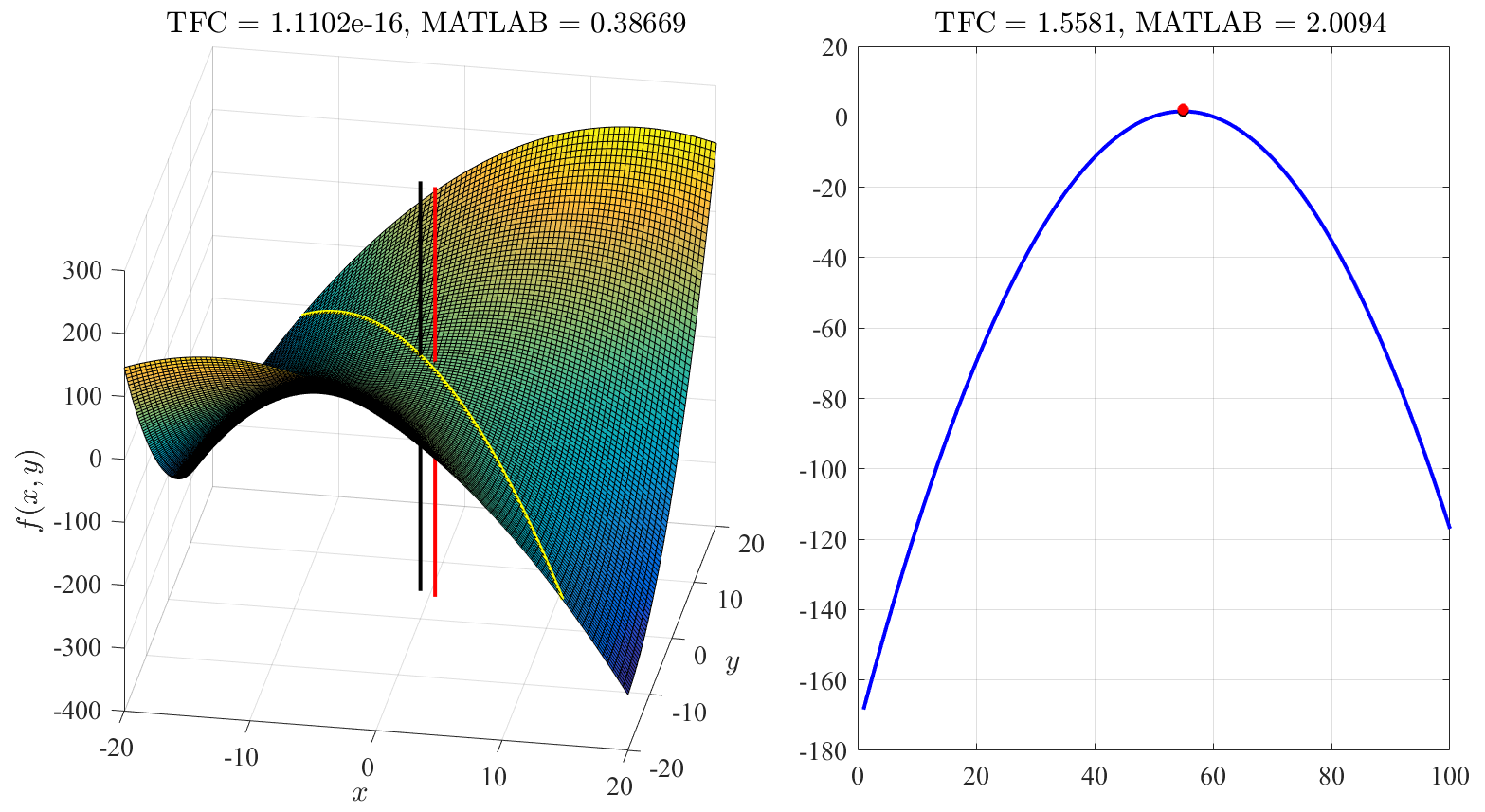}
    \caption{Monte Carlo results for QP with second constrained expression.}
    \label{fig1}
\end{figure}

Figure \ref{fig1} shows the results obtained by the TFC as null-space method (black line) and the built-in function \texttt{quadprog} of MATLAB version R2016b (red line) using random matrices $Q$ and $A$ and random vectors $\B{c}$ and $\B{b}$ for a non-convex problem (return flag $=-6$). Using random input matrices and vectors, \texttt{quadprog} may fail for several reasons: 1) Maximum number of iterations exceeded (flag $= 0$), 2) No feasible point found (flag $=-2$), 3) Problem is unbounded (flag $=-3$), 4) Non-convex problem detected (flag $=-6$, Interior-point-convex only), 5) Change in objective function is too small (flag $= 3$, Trust-region-reflective only), 6) Current search direction is not descent, no further progress can be made (flag $=-4$, Trust-region-reflective only), 7) Local minimizer found (flag $= 4$, Active-set only), and 8) Magnitude of search direction became too small, no further progress can be made. The problem is ill-posed or badly conditioned  (flag $=-7$, Active-set only).

The version R2019a of MATLAB \texttt{quadprog} is an improved version over the R2016b one. However, the algorithms adopted are still the same (MATLAB descriptions):
\begin{itemize}
  \item {\bf interior-point-convex}. This algorithm attempts to follow a path that is strictly inside the constraints. It uses a presolve module to remove redundancies and to simplify the problem by solving for components that are straightforward. The algorithm has different implementations for a sparse Hessian matrix $\B{H}$ and for a dense matrix.
  \item {\bf trust-region-reflective}. This algorithm is a subspace trust-region method based on the interior-reflective Newton's method described in \cite{quadprog}. Each iteration involves the approximate solution of a large linear system using the method of preconditioned conjugate gradients (PCG).
\end{itemize}
The current R2019a MATLAB version of \texttt{quadprog} fully satisfies the equality constraints.  Therefore, unfair\footnote{The \texttt{quadprog} is a MATLAB built-in function, which is coded in C and (most likely) highly optimized; whereas the proposed QP solution by TFC is computed by using a MATLAB script, which is line-by-line interpreted by MATLAB for its execution.} speed tests have produced the time results shown in Table \ref{speed2}, comparing the TFC with the current version of MATLAB's \texttt{quadprog}.
\begin{table}[ht]
  \centering
  \begin{tabular}{|r|r|c|c|r|}
    \hline
    $n$ & $m$ & $\Delta t$ (TFC) & $\Delta t$ (MATLAB) & time ratio \\ \hline\hline
    10 &  2 & 0.027746 & 0.87057 & 31.3759 \\
    10 &  4 & 0.029558 & 0.92745 & 31.3768 \\
    10 &  8 & 0.032127 & 0.92885 & 28.9121 \\
    20 &  4 & 0.044136 & 0.77945 & 17.6601 \\
    20 &  8 & 0.051081 & 0.94080 & 18.418 \\
    20 & 16 & 0.085394 & 0.94296 & 11.0425 \\
    40 &  8 & 0.088137 & 0.84634 & 9.6026 \\
    40 & 16 & 0.131900 & 0.81007 & 6.1415 \\
    40 & 32 & 0.198200 & 0.82851 & 4.1802 \\
    80 & 16 & 0.273250 & 0.96279 & 3.5235 \\
    80 & 32 & 0.394500 & 1.11990 & 2.8388 \\
    80 & 64 & 0.677070 & 1.32050 & 1.9502 \\
    \hline
  \end{tabular}
  \caption{Speed tests for QP with the second constrained expression.}\label{speed2}
\end{table}
The average elapsed time $\Delta t$ (in msec) in Table \ref{speed2} is obtained via $10000$ tests, where the input data ($\B{Q}$, $\B{c}$, $\B{A}$, and $\B{b}$) were randomly generated for various values of $n$ and $m$ (the number of variables and equality constraints, respectively). The time does not take into account the elapsed time to generate the input variables. Therefore, it represents only the time for the algorithm to solve QP problem. To quantify the speed gain obtained by using the TFC over the MATLAB's \texttt{quadprog}, the time ratio gain is also provided.

\vspace{10pt}

We remark here that due to the reduction in sizes of both the free vector and coefficient matrix, this second approach (with the second constrained expression) is much faster than the first approach (in Section \ref{ibelt}) while giving the same exact solution.


\section{Nonlinear Programming (NLP)}\label{nlps}

First, we introduce our framework, thanks to \cite{U-ToC, M-ToC}.  The notation will be introduced in the next section.  Given $f: \mathbb{R}^n \to \mathbb{R}$,
a nonlinear objective function that is convex and twice continuously differentiable (which also implies that \textbf{dom}$ \, f$ is open).  We will consider the constrained optimization problem
\begin{equation}\label{minf}
    \minx_{\B{x} \in \mathbb{R}^n} \;\; f (\B{x}) \qquad : \;\; \B{A} \, \B{x} \;=\; \B{b},
\end{equation}
where $\B{A} \in \mathbb{R}^{m \times n}\,, \text{ rank}(\B{A}) = m < n\,, \text{ } \B{b} \in \mathbb{R}^m$.  Denote $r = n - m$.  Let
\begin{equation}\label{x}
    \B{x} = \B{x}_0 + \B{N} \B{g},
\end{equation}
where $\B{g} \in \mathbb{R}^r\,,$ 
$\B{x}_0 \in \mathbb{R}^n$ is {\it any} particular solution of $\B{A} \, \B{x} = \B{b}\,,$ and $\B{N} \in \mathbb{R}^{n \times r}$ is a matrix having columns as a basis for the nullspace of matrix $\B{A} \in \mathbb{R}^{m \times n}\,.$ The null basis $\B{N}$ for $\B{A}$ can be computed by any of the available methods.

Note that if rank$(\B{A}) = m = n - 1$ (so $r = n - m = 1$), then the free vector $\B{g} \in \mathbb{R}^r$ becomes a free scalar variable $g \in \mathbb{R}\,.$  In this situation, the multivariate function $f (\B{x})$ can be transformed into a univariate function in $g\,;$ also, if $f(\B{x})$ is a multivariate polynomial (as a special case of functions), then it becomes a univariate algebraic expression in $g\,.$
%

\vspace{10pt}

For the equality constrained NLP problem \eqref{minf}, there are several important classical approaches \cite{Boyd}: \textit{analytically} finding a solution for its corresponding KKT system of equations, \textit{eliminating the equality constraints} to reduced the constrained problem to an unconstrained one (possessing fewer variables),
solving the \textit{dual problem} utilizing unconstrained optimization algorithms, the Newton's method with equality constraints, and combining the Newton's method with elimination scheme.  Among those, the last combination approach is in the spirit of \TFC\ (TFC) that we consider as follows.

\subsection{Second-order Newton Approach}\label{2ndnewton}

\noindent \textbf{Iteration.}

\noindent \textbf{Case} $k = 0$.  We choose $\B{g}^{(0)} \in \mathbb{R}^r$ so that $\B{N} \B{g}^{(0)} = \B{0}$.  By (\ref{x}), $\B{x}^{(0)} = \B{x}_0 + \B{N} \B{g}^{(0)}  = \B{x}_0$.

\noindent {\bf Remark.} We can take any starting point $\B{g}^{(0)} \in \mathbb{R}^r$, and the corresponding starting point is $\B{x}^{(0)} = \B{x}_0 + \B{N} \B{g}^{(0)}$.  However, for simplicity, we pick $\B{x}^{(0)} = \B{x}_0$ from \eqref{Moore-Penrose}.

Calling $\hat{f}$ the second-order Taylor approximation of $f$ at $\B{x}^{(0)}= \B{x}_0$, we have
\begin{equation*}
    f(\B{x}) = \hat{f}(\B{x}) + \text{HOT} = f(\B{x}_0) + (\B{z}_0\Ts) (\B{x} - \B{x}_0) + \dfrac{1}{2}(\B{x} - \B{x}_0)\T (\B{H}_0)  (\B{x} - \B{x}_0) + \text{HOT}\,,
\end{equation*}
where ``HOT'' stands for ``higher order terms''.
This equation, thanks to (\ref{x}), can be written as
\begin{equation}\label{t2h}
      h (\B{g}) = \hat{h} (\B{g}) + \text{HOT} = f (\B{x}_0) + \B{z}_0\Ts \, \B{N} \B{g} + \dfrac{1}{2} \, \B{g}\T \B{N}\T \, \B{H}_0 \, \B{N} \B{g} + \text{HOT}.
\end{equation}
Here, the gradient and Hessian of $f$ are respectively
\begin{align}\label{jh}
	\begin{split}
    \B{z}_k = \nabla f(\B{x}^{(k)}) = \begin{Bmatrix} \dfrac{\partial f}{\partial x_1} \\ \vdots\\ \dfrac{\partial f}{\partial x_n}\end{Bmatrix}_{\B{x}^{(k)}}\,, \quad 
    \B{H}_k = \nabla^2 f(\B{x}^{(k)}) = \begin{bmatrix} \dfrac{\partial^2 f}{\partial x_1^2} & \cdots & \dfrac{\partial^2 f}{\partial x_1 \partial x_n}\\ \vdots & \ddots & \vdots \\ \dfrac{\partial^2 f}{\partial x_n \partial x_1} & \cdots & \dfrac{\partial^2 f}{\partial x_n^2} \end{bmatrix}_{\B{x}^{(k)}}\,,
\end{split}
\end{align}
for $k = 0, 1, \cdots \,.$
Let set
\begin{equation}\label{efk}
    \begin{array}{ll}
        \B{e}_k = \left. \nabla h(\B{g}) \right|_{\B{g} = \B{g}^{(k)}} = \B{N}\T \B{z}_k, & 
        \qquad 
        \B{F}_k = \left. \nabla^2  h(\B{g}) \right|_{\B{g} = \B{g}^{(k)}} = \B{N}\T \B{H}_k \B{N} \,.
    \end{array}
\end{equation}
The first-order optimal condition for (\ref{t2h}) reads
\begin{equation*}
	\nabla \hat{h} (\B{g}) = \B{0} \qquad \Leftrightarrow \qquad \B{g}^{(1)} = -(\B{N}\T \B{H}_0 \B{N})^{-1} \B{N}\T \B{z}_0\,.
\end{equation*}
This solution $\B{g}^{(1)}$ allows us to compute a better point $\B{x}^{(1)}$ (where we expand $f (\B{x})$):
\begin{equation*}
    \B{x}^{(1)} = \B{x}_0 + \B{N} \, \B{g}^{(1)} = \B{x}_0 - \B{N} \left(\B{N}\T \, \B{H}_0 \, \B{N}\right)^{-1} \B{N}\T \B{z}_0\,.
\end{equation*}
Then the procedure is iterated as
\begin{equation}\label{xk}
    \B{x}^{(k+1)} = \B{x}_0 + \B{N} \, \B{g}^{(k+1)} =\B{x}_0 - \B{N} \sum_{j=0}^k  (\B{N}\T \B{H}_j \B{N})^{-1} \B{N}\T \B{z}_j\,.
\end{equation}

Note that it is appealing to try to simplify \eqref{xk}, but since $\B{N}$ is not a square matrix, we can not have the inverse $\B{N}^{-1}$ nor  $(\B{N}\T \B{H}_j \B{N})^{-1} = \B{N}^{-1} (\B{H}_j)^{-1} (\B{N}\T)^{-1}\,.$

\subsection{Full nonlinear Newton Approach}

Let
\begin{equation*}
	{\cal L} (\B{g}) := \nabla \hat{h} (\B{g})\,.
\end{equation*}

The stationary value can also be found using nonlinear Newton iterations
\begin{equation*}
    {\cal L} (\B{g})= \B{0}_{(n-m) \times 1}\,,
\end{equation*}
where the $(k+1)$th iteration is
\begin{equation}\label{fullnewton}
    \B{g}^{(k+1)} = \B{g}^{(k)} - (\B{F}_k)^{-1} \B{e}_k = - \sum_{j=0}^k (\B{F}_j)^{-1} \B{e}_j.
\end{equation}
Once the convergence has been obtained, that is, when $\left\|\left(\B{g}^{(k+1)} - \B{g}^{(k)}\right)\right\|_2 < \varepsilon_g$ or $\left\|{\cal L} \left(\B{g}^{(k)}\right)\right\|_2 < \varepsilon_{{\cal L}}$ (for given very small $\varepsilon_g\,, \varepsilon_{{\cal L}}$), the solution is $\B{x}^{(k+1)} = \B{x}_0 + N \, \B{g}^{(k+1)}$.  


\section{Convergence Analysis of NLP using the \TFC}\label{converge}

Applying the \TFC\ (TFC) to the constrained optimization problem (\ref{minf}), we obtain the unconstrained optimization problem
\begin{equation}\label{minh}
    \text{minimize} \quad h (\B{g}), \quad \text{for all} \quad \B{g} \in \mathbb{R}^r.
\end{equation}

\subsection{Minimization problem}

In this section, based on \cite{Boyd} (Chapter 9) by Boyd and Vandenberghe, we extend their discussion in details on the convergence analysis of the Newton's method in optimization to the minimization problem (\ref{minh}) obtained by the TFC.  We will assume that the problem (\ref{minh}) is solvable, that is, there exists an optimal point $\B{g}^*$.  (The assumptions later in this part will make sure that $\B{g}^*$ exists and is unique.)  The optimal value is denoted by $q^* = \inf\limits_{\B{g}} h(\B{g}) = h(\B{g}^*)$.

Because $h$ is convex and twice continuously differentiable, a point $\B{g}^*$ is optimal if and only if
\begin{equation}\label{critical}
    \nabla h(\B{g}^*) = 0.
\end{equation}
Therefore, solving the unconstrained minimization problem (\ref{minh}) is as finding a solution of Eq.\ (\ref{critical}).  Here, we solve the problem by an iterative algorithm, which computes a sequence of points $\B{g}^{(0)}$, $\B{g}^{(1)}$, $\cdots$ $\in$ $\textbf{dom}\, h$ so that $h\left(\B{g}^{(k)}\right) \to q^*$ as $k \to \infty$.  Such a sequence of points $\left \{\B{g}^{(k)} \right \}$ is called a \textit{minimizing sequence} for the problem (\ref{minh}).  This algorithm is ended when $h\left(\B{g}^{(k)}\right) - q^* \leq \varepsilon$, for some chosen tolerance $\varepsilon > 0$.

\noindent \textbf{Initial point and sublevel set.}

The method we are using in this part requires a proper starting point $\B{g}^{(0)}$.  That is, $\B{g}^{(0)}$ must belong to $\textbf{dom} \; h = \mathbb{R}^r$. In addition, since $h$ is continuous with $\textbf{dom}\, h = \mathbb{R}^r$ (a closed set), it follows that $h$ is closed.  Thus, by definition of closed function, for each $\B{g}^{(0)} \in \mathbb{R}^r$, the sublevel set
\begin{equation}\label{sls}
    \B{S} = \left \{ \B{g} \in \textbf{dom} \, h \; \middle| \; h(\B{g}) \leq h\left(\B{g}^{(0)}\right) \right \}
\end{equation}
is closed.  Hence, if $\textbf{dom}\, h = \mathbb{R}^r$, the initial sublevel set condition holds for any $\B{g}^{(0)} \in \textbf{dom}\, h = \mathbb{R}^r$. Equivalently, the sublevel set
\begin{equation*}
    \B{S}' = \left \{ \B{x} \in \textbf{dom} \, f \; \middle | \; f(\B{x}) \leq f \left (\B{x}^{(0)}\right), \B{A} \B{x} = \B{b} \right \}
\end{equation*}
is closed, where $\B{x}^{(0)} \in \textbf{dom} \, f$ satisfies $\B{A} \B{x}^{(0)} = \B{b}$.  This is the case if $f$ is (continuous on $\textbf{dom} \, f = \mathbb{R}^n$, thus) closed.

For $\B{x}, \B{y} \in \B{S}'$, we assume that $\nabla^2 f$ satisfies the Lipschitz condition
\begin{equation}\label{lipf}
    \left\| \nabla^2 f(\B{x}) - \nabla^2 f(\B{y})\right\|_2 \leq L \left\| \B{x} - \B{y} \right\|_2\,,
\end{equation}
where $L$ is a positive real constant.

On $\B{S}'$, we assume that $f$ is strongly convex with constant $m' > 0$, that is,
\begin{equation*}
    \nabla^2 f(\B{x}) \succeq m' \B{I}\,,
\end{equation*}
which means that $\left(\nabla^2 f(\B{x}) - m' \B{I}\right)$ is positive semi-definite (see \cite{Boyd}, p.\ 43), where $\B{I}$ is the identity matrix.

Thanks to (\cite{Boyd}, p.\ 460), the strong convexity assumption on the set $\B{S}'$ also implies that there exists $M' >0$ such that
\begin{equation}\label{uph}
    \nabla^2 f(\B{x}) \preceq M' \B{I}\,.
\end{equation}
On $\B{S}'$, we also assume that the inverse of the KTT matrix is bounded (\cite{Boyd}, p.\ 530):
\begin{equation}\label{kup}
    \left\Vert\begin{bmatrix} \nabla^2 f (\B{x}) & \B{A}\T \\ \B{A} & \B{0}_{m\times m}\end{bmatrix}^{-1}\right\Vert_2 \leq \kappa \,,
\end{equation}
where $\kappa$ is some positive constant, $\nabla^2 f (\B{x})$ is the $n\times n$ Hessian of $f (\B{x})$, and the norm $\| \cdot \|_2$ means $\sigma_{\text{max}}$ (the largest singular value of the inside matrix $\B{R}\,,$ that is, the largest square root of the eigenvalue of $\B{R}\T \B{R}$), 
which is less than or equal its Frobenius norm $\| \cdot \|_{\tu{F}}$.  Due to the symmetry of the KKT matrix, the condition \eqref{kup} indicates that its eigenvalues are bounded away from zero.

Thanks to (\cite{Boyd}, p.\ 530-531), these inequalities (\ref{uph}) and (\ref{kup}) imply that on $S$,
\begin{equation}\label{mlow}
    \nabla^2 h(\B{g}) \succeq m \B{I}\,,
\end{equation}
for some positive constant $m$.  More specifically, $m = \dfrac{\sigma_{\rm min}^2 (\B{N})}{\kappa^2 M'} $, which is positive, as $\B{N}$ is full rank. Again, thanks to (\cite{Boyd}, p.\ 460), the inequality (\ref{mlow}) on $\B{S}$ implies that there exists $M > 0$ such that
\begin{equation}\label{Mup}
    \nabla^2 h(\B{g}) \preceq M \B{I}\,.
\end{equation}

\noindent \subsubsection{Strong convexity and implications.}

\noindent  In the rest of this part, we use the assumption (\ref{mlow}) that the objective function $h(\B{g})$ is \textit{strongly convex} on $\B{S}$, that is, there exists an $m > 0$ such that
\begin{equation*}
    \nabla^2 h(\B{g}) = \B{N}\T \nabla^2 f(\B{x}) \B{N}  \succeq m \B{I} \,.
\end{equation*}
for all $\B{g} \in \B{S}$ (and with corresponding condition for $\B{x} = \B{x}_0 + \B{N} \B{g}$).  For later use, we will consider an interesting property of strong convexity.  For $\B{g}, \B{g}^{(l)} \in \B{S}$, it follows that
\begin{equation*}
    \hat{h} (\B{g}) = \hat{h} \left(\B{g}^{(l)} + \B{g} - \B{g}^{(l)}\right) = h\left(\B{g}^{(l)}\right) + \B{e}_l\T \left(\B{g}- \B{g}^{(l)}\right) + \frac{1}{2}\left(\B{g} - \B{g}^{(l)}\right)\T \B{F}_l \left(\B{g} - \B{g}^{(l)}\right).
 \end{equation*}
By the strong convexity assumption $\B{F}_l \succeq m \B{I}$, the last term on the right hand side is at least
\begin{equation*}
    \dfrac{m}{2} \left\| \B{g} - \B{g}^{(l)} \right\|^2_2.
\end{equation*}
Therefore,
\begin{equation}\label{scs}
    h(\B{g}) \geq \hat{h}(\B{g}) \geq h\left(\B{g}^{(l)}\right) + \B{e}_l\T \left(\B{g} - \B{g}^{(l)}\right) + \dfrac{m}{2} \left\| \B{g} - \B{g}^{(l)} \right\|^2_2.
\end{equation}
We thus obtain
\begin{equation*}
    \dfrac{1}{2} \left(\frac{\partial \left\|\B{g} - \B{g}^{(l)} \right\|^2_2}{\partial \B{g}} \right) = \left(\B{g} - \B{g}^{(l)}\right)\T.
\end{equation*}
For $m = 0$, (\ref{scs}) is as the basic inequality describing convexity.  When $m > 0$, we use (\ref{scs}) to bound $h\left(\B{g}^{(l)}\right) - q^*$, which is the sub-optimality of the point $\B{g}^{(l)}$, in term of $\| \B{e}_l \|_2$, where $h(\B{g}^*) = q^*$.  The right hand side of (\ref{scs}) is a convex quadratic function of $\B{g}$ (for fixed $\B{g}^{(l)}$).
Letting the gradient of this function with respect to $\B{g}$ be zero, we get that $\tilde{\B{g}} = \left( \B{g}^{(l)} - \dfrac{1}{m} \B{e}_l \right)$ minimizes the right hand side of (\ref{scs}).  Hence, Eq.\ (\ref{scs}) becomes
\begin{align*}
\begin{split}
    h (\B{g}) & \geq h \left(\B{g}^{(l)}\right) + \B{e}_l\T \left(\tilde{\B{g}} - \B{g}^{(l)}\right) + \frac{m}{2} \left\| \tilde{\B{g}} - \B{g}^{(l)} \right\|_2^2  \\ 
     & = h \left(\B{g}^{(l)}\right) + \B{e}_l\T \left (- \frac{1}{m} \B{e}_l \right) + \frac{m}{2} \left \| - \frac{1}{m} \B{e}_l \right \|_2^2 \\
    & = h\left(\B{g}^{(l)}\right) - \frac{1}{2m} \|\B{e}_l\|_2^2\,.
 \end{split}
\end{align*}
Since this holds for all $\B{g} \in \mathbb{R}^r$, so does for $\B{g} = \B{g}^*$.  Thus, we have
\begin{equation}\label{subop1}
    h \left(\B{g}^{(l)}\right) - q^* \leq  \frac{1}{2m} \left\|\B{e}_l\right\|_2^2,
\end{equation}
as desired.  This inequality demonstrates that if the gradient of $h$ is small at a point, then the point is approximately optimal.  We can also interpret the inequality (\ref{subop1}) as a condition for
\textit{sub-optimality}, which generalizes the optimality condition (\ref{critical}):
\begin{equation}\label{opt1}
    \left\| \nabla h \left(\B{g}^{(l)}\right)\right\|_2 \leq \sqrt{2 m \varepsilon} \qquad \to \qquad h\left(\B{g}^{(l)}\right) - q^* \leq \varepsilon.
\end{equation}

\noindent \textbf{The strong convexity constants}

\noindent In practice, the constants $m$ and $M$ are known only in few cases, so the inequality (\ref{opt1}) cannot be used as a termination criterion, it can only be viewed a \textit{conceptual} stopping criterion.  If we end an algorithm when $\| \nabla h\left(\B{g}^{(k)}\right) \|_2 \leq \eta$, where $\eta$ is very small, smaller than $\sqrt{m \varepsilon}$, then we obtain $h\left(\B{g}^{(k)}\right) - q^* \leq \varepsilon$.

In the following sections regarding our convergence proof for Newton's method in optimization, we will include bounds on the number of iterations needed before $h\left(\B{g}^{(k)}\right) - q^* \leq \varepsilon $, where $\varepsilon > 0$ is some tolerance.  These bounds often require the (unknown) constants $m$ and $M$, so the same remarks concern. These are at least conceptually helpful results; they deduce that the algorithm converges, even if the bound on
the number of iterations needed to reach a given exactness depends on unknown constants ($m$ or $M$).

\subsection{Descent methods}

The algorithms used in this part lead to a minimizing sequence $\B{g}^{(k)}, \ k=1, \cdots$, where
\begin{equation}\label{sd}
    \B{g}^{(k+1)} = \B{g}^{(k)} + t^{(k)} \Delta \B{g}^{(k)},
\end{equation}
and $t^{(k)} > 0$ (except when $\B{g}^{(k)}$ is optimal).  Here, $\Delta \B{g}^{(k)} \in \mathbb{R}^r$ is called the \textit{step} or \textit{search direction}.  The scalar $t^{(k)}$ is called the \textit{step size} or \textit{step length} at the $k$th iteration.  In one iteration, we can use lighter notation $\B{g}^{+} = \B{g} + t \Delta \B{g}$ for (\ref{sd}).

All the methods (including the Newton's method) in this part are \textit{descent methods}, which means that $h\left(\B{g}^{(k+1)}\right) < h\left(\B{g}^{(k)}\right)$, except when $\B{g}^{(k)}$ is optimal.  Thus, for all $k$, $\B{g}^{(k)} \in \B{S} \subset \textbf{dom} \, h$.  In a descent method, from convexity, the search direction must satisfy $\nabla h\left(\B{g}^{(k)}\right)^{\T} \Delta \B{g}^{(k)} < 0$ for all $k\,.$  Such a direction $\Delta \B{g}^{(k)}$ is called a \textit{descent direction} (for $h$, at $\B{g}^{(k)}$).

\subsubsection{Backtracking line search.}

\textbf{Algorithm.}

\noindent \textbf{given} a descent direction $\Delta \B{g}$ for $h$ at $\B{g} \in \textbf{dom} \, h$, $0 < \alpha < 0.5, 0 < \beta < 1$.

\noindent $t:=1$.

\noindent \textbf{while} $h(\B{g} + t \Delta \B{g}) > h(\B{g}) + \alpha t \ \nabla h(\B{g})^{\T} \Delta\B{g}, t:= \beta t$.

\subsection{Newton's method}

\subsubsection{The Newton step.}

For $\B{g}^{(k)} \in \textbf{dom} \, h$, the vector $\Delta \B{g}_{\tu{nt}}^{(k)} = \B{g}^{(k+1)} - \B{g}^{(k)}= - (\B{F}_k)^{-1} \B{e}_k$ is called the \textit{Newton step} or \textit{Newton direction} (for $h$, at $\B{g}^{(k)}$).  Positive definiteness of $\B{F}_k$ implies that
\begin{equation*}
    \B{e}_k^{\T} \Delta \B{g}_{\tu{nt}}^{(k)} = -\B{e}_k^{\T}(\B{F}_k)^{-1} \B{e}_k < 0
\end{equation*}
unless $\B{e}_k =0$.  Hence, the Newton step is a descent direction (unless $\B{g}^{(k)}$ is optimal).

\noindent \textbf{Minimizer of second-order approximation.}  The second-order approximation $\hat{h}$ of $h$ at $\B{g}^{(k)}$ is
\begin{equation*}
    \hat{h}\left(\B{g}^{(k)} + \B{v}\right) = h\left(\B{g}^{(k)}\right) + \B{e}_k^{\T} \B{v} + \dfrac{1}{2} \B{v}^{\T} \B{F}_k \B{v}\,,
\end{equation*}
which is a convex quadratic function of $\B{v}$ and minimized when $\B{v} = \Delta \B{g}_{\tu{nt}}^{(k)}$.

This analysis gives us some intuition about the Newton step $\Delta \B{g}_{\tu{nt}}^{(k)}$.  If the function $h$ is quadratic ($h$ = $\hat{h}$), then $\left(\B{g}^{(k)} + \Delta \B{g}_{\tu{nt}}^{(k)}\right)$ is the precise minimizer of $h$.  We will see afterward that if the function $h$ is closely quadratic, then $\left(\B{g}^{(k)} + \Delta \B{g}_{\tu{nt}}^{(k)}\right)$ should be a very good approximation of the minimizer $\B{g}^{*}$ of $h$, especially when $\B{g}^{(k)}$ is near $\B{g}^{*}$.

\noindent \textbf{Steepest descent direction in Hessian norm.}  The Newton step is also the steepest descent direction at $\B{g}^{(k)}$, for the quadratic norm defined by the Hessian $\B{F}_k$, that is,
\begin{equation}\label{hessnorm}
    \|\B{u}\|_{\B{F}_k} = \left(\B{u}\T \B{F}_k \B{u}\right)^{1/2}.
\end{equation}
This is another intuition about why the Newton step should be a good search direction, which is a very good one when $\B{g}^{(k)}$ is near $\B{g}^{*}$.

\noindent \textbf{Solution of linearized optimality condition.}  Linearizing the optimal condition $\nabla h\left(\B{g}^{*}\right) = \B{0}$ around $\B{g}^{(k)}$, we get $\nabla h\left(\B{g}^{(k)} + \B{v}\right) \approx \nabla h\left(\B{g}^{(k)}\right) + \nabla^2 h\left(\B{g}^{(k)}\right) \B{v} = \B{e}_k + \B{F}_k \B{v} = 0$, a linear equation in $\B{v}$, with solution $\B{v} = - \B{F}_k^{-1} \B{e}_k = \Delta \B{g}_{\tu{nt}}^{(k)}$.  Hence, adding the Newton step $\Delta \B{g}_{\tu{nt}}^{(k)}$ to $\B{g}^{(k)}$ makes the linearized optimality condition holds.  This again recommends that when $\B{g}^{(k)}$ is closed to $\B{g}^{*}$ (that is, the optimality conditions are almost satisfied), the update $\B{g}^{(k)} + \Delta \B{g}_{\tu{nt}}^{(k)}$ should be a very good estimate of $\B{g}^{*}$.

\noindent \textbf{Affine invariance of the Newton step.}  As a crucial character of the Newton step, the affine invariance holds.  Indeed, suppose $\B{P} \in \mathbb{R}^{n \times n}$ is nonsingular, and define $\bar{h}\left(\B{y}^{(k)}\right) = h\left(\B{P}\B{y}^{(k)}\right)$.  Then, we obtain
\begin{equation*}
    \nabla \bar{h}\left(\B{y}^{(k)}\right) = \B{P}^{\T} \B{e}_k, \qquad \nabla^2 \bar{h}\left(\B{y}^{(k)}\right) = \B{P}^{\T} \B{F}_k \B{P},
\end{equation*}
where $\B{g}^{(k)}= \B{P} \B{y}^{(k)}$.  The Newton step for $\bar{h}$ at $\B{y}^{(k)}$ is hence
\begin{align*}
    \Delta \B{y}_{\tu{nt}}^{(k)} = - \left(\B{P}^{\T} \B{F}_k \B{P}\right)^{-1} \left(\B{P}^{\T} \B{e}_k\right) =- \B{P}^{-1} \B{F}_k^{-1} \B{e}_k = \B{P}^{-1} \Delta \B{g}_{\tu{nt}}^{(k)},
\end{align*}
where $\Delta \B{g}_{\tu{nt}}^{(k)}$ is the Newton step for $h$ at $\B{g}^{(k)}$.  Therefore, the Newton steps of $h$ and $\bar{h}$ are related by the same linear transformation, and $\B{g}^{(k)} + \Delta \B{g}_{\tu{nt}}^{(k)} = \B{P} \left(\B{y}^{(k)} + \Delta \B{y}_{\tu{nt}}^{(k)}\right)\,.$

\noindent \textbf{The Newton decrement.}  We call
\begin{equation}\label{ldd}
    \lambda(\B{g}^{(k)}) = (\B{e}_k^{\T} \B{F}_k^{-1} \B{e}_k)^{1/2}
\end{equation}
the \textit{Newton decrement} (at $\B{g}^{(k)}$).  It is important to the analysis of Newton's method, and is also a stopping criterion.  The Newton decrement can connect to the difference $h\left(\B{g}^{(k)}\right) - \inf_{\B{g}} \hat{h}(\B{g})$, where $\hat{h}$ is the second order approximation of $h$ at $\B{g}^{(k)}$, as follows:
\begin{align*}
    h\left(\B{g}^{(k)}\right) - \inf_{\B{g}} \hat{h}(\B{g}) &= h\left(\B{g}^{(k)}\right) - \hat{h}\left(\B{g}^{(k)} + \Delta \B{g}_{\tu{nt}}^{(k)}\right) \\ &= - \B{e}_k^{\T} \left(\Delta \B{g}_{\tu{nt}}^{(k)}\right) - \frac{1}{2} \left(\Delta \B{g}_{\tu{nt}}^{(k)}\right)^{\T} \B{F}_k \left(\Delta \B{g}_{\tu{nt}}^{(k)}\right) \\ &= \B{e}_k^{\T} \left(\B{F}_k^{-1} \B{e}_k\right) -\frac{1}{2} \left(\B{F}_k^{-1} \B{e}_k\right)^{\T} \B{F}_k \left(\B{F}_k^{-1} \B{e}_k\right) \\ & = \frac{1}{2} \B{e}_k^{\T} \B{F}_k^{-1} \B{e}_k = \frac{1}{2} \lambda\left(\B{g}^{(k)}\right)^2.
\end{align*}
Hence, $\dfrac{1}{2} \lambda \left(\B{g}^{(k)}\right)^2$ is an approximation of $h\left(\B{g}^{(k)}\right) - q^*$, thanks to the quadratic model $\hat{h}$ of $h$ at $\B{g}^{(k)}$.  This provides a measure of the proximity of $\B{g}^{(k)}$ to $\B{g}^*$.

The Newton decrement can also be expressed as
\begin{align}\label{ldd2}
    \lambda \left(\B{g}^{(k)}\right) = \left(\left(\Delta \B{g}_{\tu{nt}}^{(k)}\right)^{\T} \B{F}_k \ \Delta \B{g}^{(k)}_{\tu{nt}}\right)^{1/2} = \left \|\Delta \B{g}^{(k)}_{\tu{nt}} \right \|_{\B{F}_k},
\end{align}
which is the Hessian norm defined by (\ref{hessnorm}).

Also, the Newton decrement occurs as a constant used in backtracking line search in the manner
\begin{equation}\label{backtrack}
    - \lambda \left (\B{g}^{(k)}\right)^2 =\B{e}_k^{\T} \left(\Delta \B{g}^{(k)}_{\tu{nt}} \right),
\end{equation}
which is the directional derivative of $h$ at $\B{g}^{(k)}$ in the direction of the Newton step:
\begin{equation*}
    - \lambda \left(\B{g}^{(k)}\right)^2 =\B{e}_k^{\T} \left(\Delta \B{g}^{(k)}_{\tu{nt}} \right) = \left. \frac{d}{\dd t} h\left(\B{g}^{(k)} + \left(\Delta \B{g}_{\tu{nt}}^{(k)}\right) t\right) \right \rvert_{t=0}.
\end{equation*}
Last, the Newton decrement is affine invariant.  That is, the Newton decrement of $\bar{h}\left(\B{y}^{(k)}\right) = h\left(\B{P}\B{y}^{(k)}\right)$ at $\B{y}^{(k)}$, where $\B{P}$ is nonsingular, is the same as the Newton decrement of $h$ at $\B{g}^{(k)} = \B{P}\B{y}^{(k)}$.

\subsubsection{Newton's method.}  In this part, Newton's method is referred to \textit{damped} Newton method or \textit{guarded} Newton method, which is different from the \textit{pure} Newton method with fixed step size $t=1$.

\noindent \textbf{\textit{Algorithm.}}

\noindent \textbf{given} a beginning point $\B{g}^{(k)} \in \textbf{dom} \ h$, tolerance $\varepsilon > 0$.

\noindent \textbf{repeat}
\begin{enumerate}
\item \textit{Calculate the Newton step and decrement.}
    \begin{equation*}
        \Delta \B{g}_{\tu{nt}}^{(k)} := -\B{F}_k^{-1}\B{e}_k \,; \quad \lambda^2 := \lambda \left(\B{g}^{(k)}\right)^2 = \B{e}_k^{\T} \B{F}_k^{-1} \B{e}_k.
    \end{equation*}
\item \textit{Termination criterion.}  \textbf{quit} if $\lambda^2 /2 \leq \varepsilon$.
\item \textit{Line search.}  Choose step size $t^{(k)}$ by backtracking line search.
\item \textit{Update.}  $\B{g}^{(k)}: = \B{g}^{(k)} + t^{(k)} \Delta \B{g}_{\tu{nt}}^{(k)}$.
\end{enumerate}
This is basically the descent method described in the previous subsection, having the Newton step as search direction.

\subsubsection{Convergence analysis.}

First, we note that the Hessian of $h$ is Lipschitz continuous on $S$ by the condition (\ref{lipf}).  Indeed, for all $\B{g}^{(k)}, \ \B{g}^{(j)}$ in $\B{S}$, we have
\begin{align*}
    \| \B{F}_k - \B{F}_j \|_2 = \left \| \B{N}^{\T} (\B{H}_k - \B{H}_j) \B{N}  \right \|_2 \leq \|\B{N}\|_2^2 \ L \ \left \|\B{x}^{(k)}-\B{x}^{(j)} \right \|_2 \leq L \ \| \B{N}\|_2^3 \ \left \|\B{g}^{(k)} - \B{g}^{(j)} \right \|_2\,.
\end{align*}
Let $K = L\|\B{N}\|^3_2$.  Thus,
\begin{equation}\label{liph}
    \|\B{F}_k - \B{F}_j \|_2 \leq K \left \| \B{g}^{(k)} - \B{g}^{(j)} \right \|_2.
\end{equation}
The coefficient $K$ can be viewed as a bound on the third derivative of $h$, and can be zero when $h$ is quadratic.  That is, $K$ measures how well $h$ can be estimated by a quadratic model.  This implies that $K$ can be important in the procedure of Newton's method.  Intuitively, Newton's method will be very effective for a function $h$ with slowly varying quadratic approximation, that is, with small $K$.

\bigskip

\noindent \textbf{Lemma.  }  \textit{The Newton's method converges quadratically (so does the sequence in (\ref{xk})), given $h: \mathbb{R}^r \to \mathbb{R}$ with corresponding assumptions in Section \ref{nlps} Nonlinear Programming, that is, $h$ is convex and twice continuously differentiable on $\textbf{dom} \, h = \mathbb{R}^r$ (which also contains the initial point $\B{g}^{(0)}$), and strongly convex on the sublevel set $\B{S}$ defined in (\ref{sls}), while the Hessian of $h$ is Lipschitz continuous on $\B{S}$.}

\begin{proof}
We first introduce the idea and outline of the convergence proof, then its details will be presented.  We will prove that there are numbers $\eta$ and $\gamma$ with $0 < \eta \leq \dfrac{m^2}{L \| \B{N} \|_2^3}$ and $\gamma > 0$ such that the following hold.
\begin{itemize}
\item If $\| \B{e}_k \|_2 \geq \eta$, then
    \begin{equation}\label{first}
        h\left(\B{g}^{(k+1)} \right) - h\left(\B{g}^{(k)}\right) \leq - \gamma\,.
    \end{equation}
\item If $\| \B{e}_k \|_2 < \eta$, then the backtracking line search chooses $t^{(k)} = 1$ and
    \begin{equation}\label{convergence}
        \frac{L \ \| \B{N}\|_2^3}{2m^2}\| \B{e}_{k+1} \|_2 \leq \left( \frac{L \ \| \B{N}\|_2^3}{2m^2} \ \| \B{e}_k \|_2\right)^2.
    \end{equation}
\end{itemize}

Let us investigate the meanings of the second condition.  Suppose that it holds for the $k$th iteration, that is, $\| \B{e}_k \|_2 < \eta$.  As $\eta \leq \dfrac{m^2}{L \| \B{N} \|_2^3}$, we deduce from (\ref{convergence}) that $\| \B{e}_{k+1} \|_2 \leq \dfrac{1}{2} \ \dfrac{L \ \| \B{N} \|_2^3}{m^2} \ \eta^2 \leq \dfrac{1}{2} \ \dfrac{1}{\eta} \ \eta^2 < \eta$.  Hence, at the $(k+1)$th iteration, the second condition also holds; and recursively, it holds for all further iterates, that is, for all $l \geq k$, we have $\|\B{e}_l \|_2 < \eta$.  Thus, for all $l \geq k$, the algorithm takes a full Newton step $t = 1$, and
\begin{equation}\label{lth}
    \frac{L \ \| \B{N}\|_2^3}{2m^2}\| \B{e}_{l+1} \|_2 \leq \left( \frac{L \ \| \B{N}\|_2^3}{2m^2} \ \| \B{e}_l \|_2\right)^2.
\end{equation}
Recursively, (\ref{lth}) implies that for all $l \geq k$,
\begin{equation*}
    \frac{L \ \| \B{N}\|_2^3}{2m^2}\| \B{e}_{l} \|_2 \leq \left( \frac{L \ \| \B{N}\|_2^3}{2m^2} \ \| \B{e}_k \|_2\right)^{2^{l-k}} \leq \left( \frac{1}{2} \right)^{2^{l-k}},
\end{equation*}
and therefore, from (\ref{subop1}), we obtain
\begin{equation*}
    h\left(\B{g}^{(l)}\right) - q^* \leq \frac{1}{2m}\|\B{e}_l\|_2^2 = \frac{2m^3}{L^2 \ \| \B{N} \|_2^6} \left(\frac{L \ \| \B{N}\|_2^3}{2m^2} \ \| \B{e}_{l} \|_2\right)^2 \leq \frac{2m^3}{L^2 \ \| \B{N} \|_2^6} \left(\frac{1}{2}\right)^{2^{l-k+1}}.
\end{equation*}
This inequality demonstrates that once the second condition (\ref{convergence}) holds, convergence is remarkably rapid, and is called \textit{quadratic convergence}.

The iterations in Newton's method inherently belong to two stages.  The first stage is referred to \textit{damped Newton phase} because a step size $t < 1$ can be chosen.  The second stage, which arises when the condition $\| \B{e}_k \|_2 \leq \eta$ is satisfied, is called the \textit{quadratically convergent stage} or \textit{pure Newton phase}, as a step size $t=1$ is selected.  

We now approximate the total complexity.  The number of iterations until $h\left(\B{g}^{(k)}\right) - q^* \leq \varepsilon$ has an upper bound
\begin{equation}\label{ubn}
    \frac{h\left(\B{g}^{(0)}\right) - q^*}{\gamma} + \tu{log}_2 \ \tu{log}_2 (\varepsilon_0 / \varepsilon) \approx \frac{h\left(\B{g}^{(0)}\right) - q^*}{\gamma} + 6.
\end{equation}
The first term of (\ref{ubn}) corresponds to the upper bound on the number of iterations in the damped Newton phase.  The second term of (\ref{ubn}) corresponds to the number of iterations in the quadratically convergent phase, which grows uncommonly slowly with required exactness $\varepsilon$, and can be viewed as a constant, say six, so $\varepsilon \approx 5 \cdot 10^{-20} \ \varepsilon_0$, where $\varepsilon_0 = \dfrac{2m^3}{L^2 \ \| \B{N} \|_2^6}\,.$

\noindent \textbf{Damped Newton phase.}

We now find $\gamma$ in the inequality (\ref{first}).  Assume $\| \B{e}_k \|_2 \geq \eta$.  We first search a lower bound on the step size chosen from the backtracking line search.  Strong convexity (\ref{Mup}) gives $\B{F}_k \succeq \eta$ on $\B{S}\,.$  Hence,
\begin{align*}
    h \left (\B{g}^{(k)} + t \Delta \B{g}_{\tu{nt}}^{(k)} \right) & \leq h \left (\B{g}^{(k)}\right) + t \B{e}_k^{\T} \ \Delta \B{g}_{\tu{nt}}^{(k)} + \frac{M}{2} \ \left \|\Delta \B{g}_{\tu{nt}}^{(k)} \right \|^2_2 \ t^2 \\ & \leq h \left (\B{g}^{(k)}\right) - t \lambda \left (\B{g}^{(k)} \right)^2 + \frac{M}{2m} \ t^2 \ \lambda \left (\B{g}^{(k)} \right)^2,
\end{align*}
where (\ref{backtrack}) is applied, and by (\ref{ldd2}), we have $\lambda \left (\B{g}^{(k)} \right)^2 = \left (\Delta \B{g}_{\tu{nt}}^{(k)} \right)^{\T} \ \B{F}_k \ \Delta \B{g}^{(k)}_{\tu{nt}} \geq m \ \left \| \Delta \B{g}^{(k)}_{\tu{nt}} \right \|_2^2$.  Since $\|\B{F}_k^{-1} \|_2 \geq \| \B{F}_k \|_2^{-1}$, we note also that
\begin{equation}\label{Mieq}
 \lambda \left (\B{g}^{(k)} \right)^2 = \B{e}_k^{\T} \B{F}_k^{-1} \B{e}_k \geq \frac{1}{M} \| \B{e}_k \|_2^2\,.
\end{equation}

The step size $\hat{t} = m/M$ satisfies the exit condition of the line search, as
\begin{align*}
    h \left (\B{g}^{(k)} + \hat{t} \Delta \B{g}_{\tu{nt}}^{(k)} \right)  \leq h \left (\B{g}^{(k)} \right) - \frac{m}{2M} \lambda \left (\B{g}^{(k)} \right)^2 \leq h \left (\B{g}^{(k)}\right) - \alpha \ \hat{t} \ \lambda \left (\B{g}^{(k)} \right)^2\,,
\end{align*}
where $0 < \alpha < 0.5$.
Thus, the line search yields a step size $t \geq \beta \ m/M$, $0 < \beta < 1$, leading to a decrease of the objective function
\begin{align*}
    h \left (\B{g}^{(k)} + t \Delta \B{g}_{\tu{nt}}^{(k)} \right) - h \left(\B{g}^{(k)} \right)  \leq - \alpha \ t \ \left(\lambda \left(\B{g}^{(k)}\right) \right)^2 \leq - \alpha \ \beta \ \frac{m}{M^2} \ \|\B{e}_k \|_2^2 \leq - \alpha \ \beta \ \eta^2 \ \frac{m}{M^2}\,,
\end{align*}
where we make use of (\ref{Mieq}). Hence, (\ref{first}) holds for
\begin{equation*}
    \gamma = \alpha \ \beta \ \eta^2 \ \frac{m}{M^2}\,.
\end{equation*}


\noindent \textbf{Quadratically convergent phase.}

We now derive the inequality (\ref{convergence}).  Assume $\| \B{e}_k \|_2 < \eta$.  First, we prove that the backtracking line search choose unit steps $t=1$, given $\eta \leq 3(1 - 2 \alpha) \ \dfrac{m^2}{L \ \| \B{N} \|_2^3}$.  From the Lipschitz condition (\ref{liph}), we obtain, for $t \geq 0$,
\begin{align*}
    \left \| \nabla^2 h\left( \B{g}^{(k)} + t \ \Delta \B{g}^{(k)}_{\tu{nt}}\right) - \nabla^2 h\left(\B{g}^{(k)}\right) \right \|_2 
    \leq t \ K \left \| \Delta \B{g}^{(k)}_{\tu{nt}} \right \|_2.
\end{align*}
Hence,
\begin{align*}
    \left | \left(\Delta \B{g}^{(k)}_{\tu{nt}}\right)^{\T} \ \left( \nabla^2 h\left( \B{g}^{(k)} + t \ \Delta \B{g}^{(k)}_{\tu{nt}}\right) - \nabla^2 h\left(\B{g}^{(k)} \right)\right) \Delta \B{g}^{(k)}_{\tu{nt}} \right | 
    \leq t \ K \left \| \Delta \B{g}^{(k)}_{\tu{nt}} \right \|_2^3.
\end{align*}
Let
\begin{equation}\label{htil}
    \tilde{h}(t) = h\left( \B{g}^{(k)} + t \ \Delta \B{g}^{(k)}_{\tu{nt}}\right),
\end{equation}
we arrive at
\begin{equation*}
    \tilde{h}''(t) = \left(\Delta \B{g}^{(k)}_{\tu{nt}}\right)^{\T} \  \nabla^2 h\left( \B{g}^{(k)} + t \ \Delta \B{g}^{(k)}_{\tu{nt}}\right) \ \Delta \B{g}^{(k)}_{\tu{nt}},
\end{equation*}
so the previous inequality becomes
\begin{equation}\label{preineq}
    \left |\tilde{h}''(t) - \tilde{h}''(0) \right | \leq t \ K \ \left \| \Delta \B{g}^{(k)}_{\tu{nt}} \right \|_2^3.
\end{equation}
This inequality will be used to identify an upper bound on $\tilde{h}''(t)$.  Noting that
\begin{equation*}
    \tilde{h}''(0) = \lambda \left(\B{g}^{(k)}\right)^2 = \left(\Delta \B{g}^{(k)}_{\tu{nt}}\right)^{\T} \  \nabla^2 h\left( \B{g}^{(k)} \right) \ \Delta \B{g}^{(k)}_{\tu{nt}}  \geq m \ \left \| \Delta \B{g}^{(k)}_{\tu{nt}} \right \|_2^2,
\end{equation*}
we begin with the following inequality (derived from the inequality (\ref{preineq})):
\begin{align*}
    \tilde{h}''(t) \leq \tilde{h}''(0) + t \ K \ \left \| \Delta \B{g}^{(k)}_{\tu{nt}} \right \|_2^3 \leq \lambda \left(\B{g}^{(k)}\right)^2 + t \ \dfrac{K}{m^{3/2}} \lambda \left (\B{g}^{(k)}\right)^3.
\end{align*}
Integrating this inequality, we obtain
\begin{align*}
    \tilde{h}'(t) & \leq \tilde{h}(0) + t \ \lambda \left (\B{g}^{(k)}\right)^2 + t^2 \ \dfrac{K}{2 \ m^{3/2}} \lambda \left(\B{g}^{(k)}\right)^3\\
    & =-\lambda \left (\B{g}^{(k)} \right )^2 + t \ \lambda \left(\B{g}^{(k)}\right)^2 + t^2 \ \dfrac{K}{2 \ m^{3/2}} \lambda \left(\B{g}^{(k)}\right)^3,
\end{align*}
where $\tilde{h}'(0) = \B{e}_k^{\T} \ \Delta \B{g}^{(k)}_{\tu{nt}} = - \B{e}_k^{\T} \ \B{F}_k^{-1} \ \B{e}_k = - \lambda \left(\B{g}^{(k)}\right)^2$ by (\ref{ldd}).  Integrating the previous inequality once more, we get
\begin{equation*}
    \tilde{h}(t) \leq \tilde{h}(0) - t\lambda \left(\B{g}^{(k)}\right)^2 + t^2 \ \dfrac{1}{2} \ \lambda \left(\B{g}^{(k)}\right)^2 + t^3 \ \dfrac{K}{6 \ m^{3/2}} \lambda \left(\B{g}^{(k)}\right)^3.
\end{equation*}
Last, picking $t=1$, from the definition (\ref{htil}) of $\tilde{h}(t)$, we get
\begin{align}\label{last}
    h \left(\B{g}^{(k)} + \ \Delta \B{g}^{(k)}_{\tu{nt}}\right) \leq  h\left(\B{g}^{(k)}\right) - \frac{1}{2} \ \lambda \left(\B{g}^{(k)}\right)^2 + \dfrac{K}{6 \ m^{3/2}} \lambda \left(\B{g}^{(k)}\right)^3\,.
\end{align}
Now, suppose $\| \B{e}_k \|_2 \leq \eta \leq 3(1 - 2\alpha) \dfrac{m^2}{L \ \|\B{N}\|_2^3}\,.$ By the strong convexity of $h$, we get
\begin{equation}\label{Finvert}
    \left \|\B{F}_k^{-1} \right \|_2 \leq \dfrac{1}{m}\,.
\end{equation}
We then obtain
\begin{align*}
    \lambda \left (\B{g}^{(k)}\right ) = \left(\B{e}_k^{\T} \ \B{F}_k^{-1} \ \B{e}_k\right)^{1/2} \leq \| \B{e}_k \|_2 \left \| \B{F}_k^{-1} \right \|_2^{1/2} \leq 3(1 - 2\alpha) \dfrac{m^{3/2}}{L \|\B{N}\|_2^3}\,.
\end{align*}
Using (\ref{last}), we get
\begin{align*}
    h\left ( \B{g}^{(k)} + \ \Delta \B{g}^{(k)}_{\tu{nt}}\right ) \leq  h\left(\B{g}^{(k)}\right) - \lambda \left(\B{g}^{(k)}\right)^2 \left(\frac{1}{2} - \dfrac{K \ \lambda \left(\B{g}^{(k)}\right)}{6 \ m^{3/2}}\right) \leq h\left(\B{g}^{(k)}\right) + \alpha \ \B{e}_k^{\T} \ \Delta \B{g}^{(k)}_{\tu{nt}},
\end{align*}
which demonstrates that the unit step $t=1$ is accepted by the backtracking line search.

Let us now investigate the rate of convergence.  Making use of the Lipschitz condition, we obtain
\begin{align*}
    \| \B{e}_{k+1} \|_2 & = \left \| \nabla h\left( \B{g}^{(k)} + \ \Delta \B{g}^{(k)}_{\tu{nt}}\right)\right \|_2\\ 
    & = \left \| \nabla h\left( \B{g}^{(k)} + \ \Delta \B{g}^{(k)}_{\tu{nt}}\right) - \nabla h\left(\B{g}^{(k)}\right) - \B{F}_k  \left( - \B{F}_k^{-1} \ \B{e}_k \right) \right \|_2 \\ 
    & = \left \| \int_0^1 \ \left(\nabla^2 h\left(\B{g}^{(k)} + t \ \ \Delta \B{g}^{(k)}_{\tu{nt}}\right) - \B{F}_k\right) \ \Delta \B{g}^{(k)}_{\tu{nt}} \ \dd t \right \|_2 \\ 
    & \leq \int_0^1 \left \| \nabla^2 h\left( \B{g}^{(k)} + t \  \ \Delta \B{g}^{(k)}_{\tu{nt}}\right) - \nabla^2 h\left(\B{g}^{(k)}\right) \right \|_2 \ \left \|\Delta \B{g}^{(k)}_{\tu{nt}}\right \|_2  \ \dd t \\ 
    & \leq \int_0^1 K \ \left \|\Delta \B{g}^{(k)}_{\tu{nt}} \right \|_2^2 \ t \ \dd t \\
    &= \frac{K}{2} \ \| \B{F}_k^{-1} \B{e}_k\|_2^2\\ &\leq \frac{L \ \| \B{N} \|_2^3}{2m^2}\ \|\B{e}_k \|_2^2\,,
\end{align*}
which is (\ref{convergence}) as desired (where the last inequality follows from (\ref{Finvert})).

In conclusion, the procedure chooses unit steps $t =1$ and the condition (\ref{convergence}) holds if $\| \B{e}_k \|_2 < \eta$, where
\begin{equation*}
    \eta = {\rm min} \left \{1, 3(1-2 \alpha)\right \} \ \frac{m^2}{L \ \| \B{N} \|_2^3}\,.
\end{equation*}
Substituting this bound and $\gamma = \alpha \ \beta \ \eta^2 \dfrac{m}{M^2}$ into (\ref{ubn}), we deduce that the number of iterations $D_{\max}$ is bounded above by
\begin{equation}\label{iterations}
    D_{\max} \leq 6 + \frac{M^2 \ L^2 \ \| \B{N} \|_2^6}{ \alpha \beta \ m^5 \ {\rm min} \left\{ 1, \ 9 (1 - 2 \alpha)^2 \right \}} \ \left(h\left(\B{g}^{(0)}\right) - q^*\right).
\end{equation}
\end{proof}

\section{Numerical example of nonlinear programming}\label{enlp} 
Consider the minimization problem \eqref{minf}
	\begin{equation}\label{minfe}
		\minx_{\B{x} \in \mathbb{R}^n} \;\; f (\B{x}) \qquad : \;\; \B{A} \, \B{x} \;=\; \B{b},
	\end{equation}
	with the following convex nonlinear function in variable $\B{x} = (x, y, z)$ as

\begin{equation}\label{egn1}
	f (\B{x}) = x y (x y + 6 y - 8 x - 48) + z^2 - 8 z + 9 y^2 - 72 y + 16 x^2 + 96 x + 160\,,
\end{equation}
subject to the specific equality constraint
\begin{equation}
	\B{A} \, \B{x} = \begin{bmatrix*}[r] 1 & 2 & -1 \\ 1 & 0 & 1\end{bmatrix*} \begin{bmatrix} x \\ y \\ z\end{bmatrix} = \begin{bmatrix*}[r] 1 \\ 1\end{bmatrix*} = \B{b}\,.
\end{equation}
Benefiting the notation \eqref{jh} from Section \ref{nlps}, the gradient of $f$ is
\begin{equation*}
	\B{z}_k = 2\begin{bmatrix} x y^2 + 3 y^{2} - 8 x y -24 y + 16 x + 48 \\ x^2 y + 6 x y + 9 y - 4 x^2 - 24 x - 36 \\ z - 4\end{bmatrix}\,,
\end{equation*}
and the Hessian of $f$ is
\begin{equation*}
	\B{H}_k = 2 \begin{bmatrix} y^2 - 8 y + 16 & 2 x y + 6 y - 8 x -24 & 0 \\ 2 x y + 6 y - 8 x - 24 & x^2 + 6 x + 9 & 0 \\ 0 & 0 & 1\end{bmatrix}\,.
\end{equation*}
Starting with $\B{g}^{(0)} = \B{0}$, that is, with $\B{x}^{(0)} = \B{x}_0 = \B{A}^{+} \,\B{b} = \begin{Bmatrix} 1, & 0, & 0\end{Bmatrix}\T$ from \eqref{Moore-Penrose}, the solution $\bfa{x}$ of \eqref{minfe} in the form \eqref{xk} can be obtained in approximately 8 iterations. The results are provided in Table \ref{tab:1}.
\begin{table}[ht]
	\centering
	\begin{tabular}{|c|c|c|c|c|c|} \hline
		Iteration $k$ & $x$ & $y$ & $z$ & $f(x, y, z)$ & $||\B{A} \, \B{x}_k - \B{b}||_2$ \\ \hline \hline
		$0$ & $1$ & $0$ & $0$ & $272$ & $0$ \\
		$1$ & $-0.36082$ & $1.3608$ & $1.3608$ & $55.4799$ &  $4.9651 \cdot 10^{-16}$ \\
		$2$ & $-1.2817$ & $2.2817$ & $2.2817$ & $11.6709$ &  $8.8818 \cdot 10^{-16}$ \\
		$3$ & $-1.9157$ & $2.9157$ & $2.9157$ & $2.5583$ &  $1.4043 \cdot 10^{-15}$ \\
		$4$ & $-2.3668$ & $3.3668$ & $3.3668$ & $0.56159$ &  $1.7764 \cdot 10^{-15}$ \\
		$5$ & $-2.7019$ & $3.7019$ & $3.7019$ & $0.096793$ &  $1.7764 \cdot 10^{-15}$ \\
		$6$ & $-2.9309$ & $3.9309$ & $3.9309$ & $0.0048027$ &  $1.7764 \cdot 10^{-15}$ \\
		$7$ & $-2.9987$ & $3.9987$ & $3.9987$ & $1.6513 \cdot 10^{-6}$ &  $1.831 \cdot 10^{-15}$ \\
		$8$ & $-3$ & $4$ & $4$ & $5.6843 \cdot 10^{-14}$ &  $1.7764 \cdot 10^{-15}$ \\
		\hline
	\end{tabular}
	\caption{Iteration results for NLP.}
	\label{tab:1}
\end{table}

\section{Discussion}\label{discuss}

Given a convex and twice continuously differentiable objective function $h\,,$ one can solve the unconstrained minimization problem \eqref{minh} by finding a solution of the corresponding optimality equation \eqref{critical}:
\begin{equation*}
	\nabla h(\B{g}) = 0\,,
\end{equation*} 
which is a set of $n$ equations in the $n$ variables $g_1, g_2, \cdots, g_n\,.$  In a few special situations, solution for this equation can be found \textit{analytically}.  Usually, however, to solve this problem, one needs an \textit{iterative} algorithm, such as the Newton's method (Sections \ref{nlps} and \ref{converge}) in the spirit of \TFC\ (TFC).  

By definition, an iterative technique produces a series of approximations to the solution of differential equations, where the current approximation is obtained from one or more prior approximations.  In this section, we will discuss the ideas how to solve the nonlinear optimality equation \eqref{critical} by other iterative methods (rather than the Newton's method), such as Picard \cite{modifiedPicard}, mixed Picard-Newton \cite{mixpn94}, and gradient methods accelerated by Picard-Mann hybrid iterative process \cite{psurvey}.  

\subsection{Picard method}\label{pd}

\noindent  We first recall the classical definition of the \textbf{Fixed Point Iteration Method}: In this approach, Eq.\ \eqref{critical} is rewritten in the form
\begin{equation}\label{fixpt}
\bfa{g} = \B{r}(\bfa{g})	
\end{equation}
such that any solution of \eqref{fixpt}, which is a fixed point of the vector-valued function $\B{r}\,,$ is a solution of \eqref{critical}.  

Now, starting from an initial guess $\B{g}^{(0)}$ (see Subsection \ref{2ndnewton}, for instance), then with $k=0, 1, \cdots\,,$ the iterative process
\begin{equation}\label{picard1}
 \B{g}^{(k+1)} = \B{r}(\B{g}^{(k)})
\end{equation}
 is called \textbf{Picard iteration}.  Note that if we take $\B{r}(\B{g}^{(k)})$ to be the right hand side of the \textbf{Newton's iteration} \eqref{fullnewton}, then the Newton's method is a particular case of Picard algorithm.  Conversely, the Picard process may be rewritten as the Newton scheme \cite{pnye21, mixpn94}.  

This Picard procedure can be used to solve the nonlinear optimal equation \eqref{critical} as for Richards equation \cite{modifiedPicard, gne, cemnlporo, Spiridonov2019, rtt21}.



\subsection{Mixed Picard-Newton method}\label{mpn}
When the Picard scheme can be written in the form of Newton method, one can benefit the best sides of them in a mixed Picard-Newton method \cite{mixpn94} .  
It was based on the idea
of utilizing Picard iteration to improve the ``initial'' solution estimate for the Newton
method. The Picard scheme is employed for the first few iterations until it begins to
converge steadily, and then the Newton scheme is used for the subsequent iterations in
this strategy. Once Picard iteration begins to converge in the implementation, the
transition to Newton iteration is completed after a certain reduction in convergence error. 

\subsection{Gradient methods accelerated by Picard-Mann hybrid iterative process}\label{pm}
This subsection is based on \cite{psurvey, hsm18} for a generally unconstrained nonlinear optimization problem, whose objective multivariate function $h$ is now assumed to be twice continuously differentiable and uniformly convex (which is thus strictly convex and then convex). 

We first discuss about the Picard-Mann hybrid iterative process \cite{pm13, psurvey, hsm18}, for finding a fixed point of the continuous function $\B{r}$ defined in \eqref{fixpt}.  Choosing an initial guess $\B{g}^{(0)}$ (see Subsection \ref{2ndnewton}, for instance), the hybrid Picard-Mann iterations \cite{psurvey} induce two following sequences $\B{g}^{(k)}\,, \B{w}^{(k)}$ by the relations (for $k = 0, 1, \cdots$):
\begin{align}\label{pmrules}
\begin{split}
\begin{cases}
\B{g}^{(k+1)} = \B{r} (\B{w}^{(k)})\,,\\
	\B{w}^{(k)} = (1- \nu_k) \B{g}^{(k)} + \nu_k \B{r} (\B{g}^{(k)}) \,,	
\end{cases}
\end{split}
\end{align}
equivalently,
\begin{align}\label{pmrules2}
\B{g}^{(k+1)} = \B{r} [(1 - \nu_k) \B{g}^{(k)} + \nu_k \B{r} (\B{g}^{(k)})] \,.	
\end{align}
Here, the real positive number $\nu_k \in (0,1)$ is called the \textit{correction parameter} in \cite{hsm18}.

In numerical studies, \cite{pm13} suggested a set of constant values $\nu = \nu_k \in (0,1)\, \forall \, k$ and found that the process \eqref{pmrules} converges quicker than the classical Picard, Mann, and Ishikawa (see the references in \cite{pm13}, for instance).  Therefore, an application of this Picard-Mann hybrid iterative process is to accelerate iterations for solving the optimal equation \eqref{critical} or more general nonlinear optimization problems \cite{psurvey}.   

In addition, to accelerate \textit{gradient methods} \cite{psurvey}, the process \eqref{pmrules} was used in \cite{hsm18} for a hybridization HSM of the SM approach \cite{sm10}.  More specifically, first, the SM algorithm \cite{sm10} is based on a transformation of the Newton with line search into a \textit{gradient descent method}, which has the features of quasi-Newton and modified Newton methods \cite{sm10, hsm18}.  Thus, the estimation to the inverse Hessian matrix in \eqref{efk} of the objective function $h$ from \eqref{minh} in SM iteration is provided by the scalar diagonal matrix $\gamma_k^{-1} \B{I}\,,$ where $\gamma_k$ is called the \textit{acceleration parameter} and is obtained by applying the Taylor expansion on the SM algorithm \cite{hsm18}.  With $h$ in \eqref{minh}, recall from \eqref{efk} that
\begin{equation*}
	\begin{array}{ll}
		\B{e}_k = \left. \nabla h(\B{g}) \right|_{\B{g} = \B{g}^{(k)}}\,,
	\end{array}
\end{equation*}
for $k = 0, 1, \cdots \,.$  Now, from \eqref{picard1}, we assume that $\B{r}(\B{g}^{(k)})$ is defined by SM iteration
\begin{equation}\label{smi}
\B{r}(\B{g}^{(k)}) = \B{g}^{(k)} - \gamma_k^{-1} t^{(k)} \B{e}_k \,,
\end{equation}
where the scalar $t^{(k)}$ is called the \textit{step size} or \textit{step length} at the $k$th iteration, as in \eqref{sd}.  Using this mapping $\B{r}$ from \eqref{smi}, the iterations \eqref{pmrules} or \eqref{pmrules2} turn into the HSM process (as hybridization of the SM method) \cite{hsm18, psurvey}:
\begin{equation}\label{hsme}
\B{g}^{(k+1)} = \B{g}^{(k)} - (\nu_k +1) \gamma_k^{-1} t^{(k)} \B{e}_k \,,
\end{equation}
where the $(k+1)$-th acceleration parameter value is defined by
\[\gamma_{k+1} = 2\gamma_{k} \frac{\gamma_{k} [h(\B{g}^{(k+1)}) - h(\B{g}^{(k)})] + (\nu_{k} + 1) t^{(k)} \| \B{e}_{k} \|^2}{(\nu_{k} + 1)^2 (t^{(k)})^2 \| \B{e}_{k} \|^2}\,.
\]
In summary, the HSM \cite{hsm18} is a gradient descent technique with a line search strategy for addressing unconstrained optimization problems.  More specifically, the HSM \cite{hsm18} is defined as a result of using a Picard-Mann hybrid iterative process on the accelerated gradient descent SM method \cite{sm10}. 

In \cite{mhsm18}, a suitable starting value was proposed in the backtracking procedure to define a modified HSM (MHSM) technique.  A hybridization of the ADD approach was studied in \cite{hadd18}.  The hybridization HTADSS has recently been proposed, researched, and tested in \cite{htadss19}. 

All the methods here (SM \cite{sm10}, HSM \cite{hsm18}, MHSM \cite{mhsm18}, HADD \cite{hadd18}, HTADSS\cite{htadss19}) were proved (in those references) that they are at least linearly convergent on the sets of uniformly convex functions and strictly convex quadratic functions.  In addition, the hybridization HSM converges faster than the initial SM \cite{hsm18}.

\section{Conclusions}\label{conclude}

This paper demonstrates how to use the \TFC\ (TFC) to solve quadratic and nonlinear programming problems subject to linear equality constraints. This is handled without employing the traditional Lagrange multiplier method.  

In particular, for equality constrained \textit{quadratic programming} (QP) problem, two efficient TFC approaches are introduced. The \textit{first TFC approach} specifies the constrained expression \eqref{eta} as in Section \ref{2belts}, that is \eqref{DM2}.  This transforms the initial constrained optimization problem into an unconstrained linear system of equations, where the variable is the free vector $\B{g}$. 
The solution is obtained by solving this system using the Moore--Penrose pseudo-inverse.  The \textit{second TFC approach} takes advantage of the classical null-space methods to express the equality constraints. In this case, the free vector $\B{g}$'s dimension is $n-m\,,$ which is the difference between $n$ (the size of the solution vector $\B{x}$) and $m$ (the rank of equality constrained matrix).  This approach also leads the constrained QP problem to an unconstrained linear system, but its coefficient matrix is now nonsingular. Because of the reduction in sizes of both the free vector and the matrix inverse, this second approach is consequently faster than the first one while giving the same exact solution.  Regarding the \textit{second approach}, numerical validation tests 
have been performed to quantify the \textit{accuracy} with respect to the TFC and R2016b MATLAB version of \texttt{quadprog}.  For each approach, we provide $10000$ \textit{speed} (Monte Carlo) tests with respect to the TFC and the improved R2019a MATLAB version of \texttt{quadprog} for solving random QP problems.

The solution of \textit{nonlinear programming} (NLP) problem subject to linear equality constraints is obtained by Newton's method combining with elimination scheme in the spirit of TFC.  The starting point is a feasible initial vector $\B{x}^{(0)}$, which is a closed-form solution of the linear equality constraints.   The nonlinear objective function is expanded into second order.  We provide convergence analysis of the proposed approach to solve the NLP problem using the TFC.  This analysis gives the termination criterion, quadratic convergence rate (\ref{convergence}), and an upper bound (\ref{iterations}) on the total number of iterations required to attain a given accuracy.  A numerical validation test has been presented for this NLP problem.


\vspace{20pt}

\noindent \textbf{Acknowledgements.}  

Tina Mai's research is funded by Vietnam National Foundation for Science and 
Technology Development (NAFOSTED) under grant number 101.99-2019.326.  The authors thank the Reviewers very much for their thoughtful comments, which have led to an improvement in the
quality of the paper.



\bigskip

\bibliographystyle{plain}
\bibliography{ReferencesTFC}

\end{document}